\documentclass[12pt]{article}

\usepackage{graphicx}
\usepackage{lscape}
\usepackage{color}
\usepackage{theorem}
\theoremstyle{plain}
 \newtheorem{theorem}{Theorem}[section]
 \newtheorem{proposition}[theorem]{Proposition}
 
 \newtheorem{lemma}[theorem]{Lemma}
 \newtheorem{corollary}[theorem]{Corollary}
 \newtheorem{example}[theorem]{Example}
 \newtheorem{remark}[theorem]{Remark}
 
 \newtheorem{definition}[theorem]{Definition}

\usepackage{amsmath, amssymb}  
\numberwithin{equation}{section}

\def\om{\omega}
\def\lam{\lambda} 
\def\threedot#1{\stackrel{...}{#1}}
\def\fourdot#1{\stackrel{....}{#1}}
\def\tom{\tilde{\om}}
\def\te{\tilde{e}}
\def\comment#1{}

\allowdisplaybreaks[1]

\title{General-affine invariants of plane curves and space curves}

\makeatletter
\renewcommand\@date{{%
  \vspace{-\baselineskip}%
  \large\centering
  \begin{tabular}{@{}c@{}}
    Shimpei Kobayashi\footnote{The first named author is partially supported by Kakenhi 26400059 and Deutsche Forschungsgemeinschaft-Collaborative Research Center, TRR 109, ``Discretization in Geometry and Dynamics''.} \textsuperscript{1} \\
    \normalsize \texttt{shimpei@math.sci.hokudai.ac.jp}
  \end{tabular}%
  \quad and\quad
  \begin{tabular}{@{}c@{}}
    Takeshi Sasaki \textsuperscript{2} \\
    \normalsize \texttt{sasaki@math.kobe-u.ac.jp}
  \end{tabular}

  \bigskip

 \textsuperscript{1}Department of Mathematics, Hokkaido University,  
 Sapporo, 060-0810, Japan 
\par
 \textsuperscript{2}Department of Mathematics, Kobe University, 
 Kobe, 657-8501, Japan

\bigskip

\today

}}
\makeatother

\makeatletter
\newcommand{\subjclass}[2][2010]{%
  \let\@oldtitle\@title%
  \gdef\@title{\@oldtitle\footnotetext{#1 \emph{Mathematics subject classification.} #2}}%
}
\newcommand{\keywords}[1]{%
  \let\@@oldtitle\@title%
  \gdef\@title{\@@oldtitle\footnotetext{\emph{Keywords.} #1.}}%
}
\makeatother

\begin{document}

\maketitle
\begin{abstract}
We present a fundamental theory of curves in the affine plane 
and the affine space, equipped with the general-affine 
groups ${\rm GA}(2)={\rm GL}(2,{\bf R})\ltimes {\bf R}^2$ 
and ${\rm GA}(3)={\rm GL}(3,{\bf R})\ltimes {\bf R}^3$, respectively.
We define general-affine length parameter and curvatures and show 
how such invariants determine the curve up to general-affine motions.
 We then study  the extremal problem of the general-affine length functional 
and derive a variational formula. We give several examples of curves and
also discuss some relations with equiaffine treatment 
and projective treatment 
of curves.
\end{abstract}

\section{Introduction}
Let ${\bf A}^n$ be the affine $n$-space with coordinates $x=(x^1,\dots,x^n)$.
It is called a unimodular affine space if it is equipped with a
parallel volume element, namely a determinant function.
The unimodular affine group
${\rm SA}(n)={\rm SL}(n,{\bf R})\ltimes {\bf R}^n$ acts as
\[ 
x=(x^i) \longrightarrow gx+a={\Big(}\sum_{j} g_j^ix^j + a^i{\Big)},
\quad g=(g_j^i)\in {\rm SL}(n,{\bf R}),\ 
a=(a^i)\in {\bf R}^n,
\]
which preserves the volume element. 
Study of geometric properties of submanifolds 
in ${\bf A}^n$ invariant under this group is called equiaffine differential 
geometry, while the study of properties 
invariant under the general affine group
${\rm GA}(n)={\rm GL}(n,{\bf R})\ltimes {\bf R}^n$ is called 
general-affine differential geometry. Furthermore, 
the study of the geometric properties of submanifolds 
in the projective space ${\bf P}^n$
invariant under the projective linear group 
${\rm PGL}(n)={\rm GL}(n+1,{\bf R})/{\rm center}$ 
is called projective differential geometry. 
Equiaffine differential geometry, as well as
projective differential geometry, has long been studied and
has yielded a plentiful amount of results,  especially for
curves and hypersurfaces; we refer to \cite{Bl,Sc,NS} and 
\cite{Wi, La, Bol} to name a few references. However, the study of
general-affine differential geometry is little known even for
curves. The purpose of this paper is to present a basic study of
plane curves and space curves in general-affine differential geometry
by recalling old results and by adding some new results.
In addition, we relate them with
the curve theory in equiaffine and projective differential geometry.
Although study of invariants of curves of higher-codimension could possibly be
given by a similar formulation used in this paper, it probably requires
a more complicated presentation and is not attempted here.
For the studies of submanifolds in the affine space
which correspond to other types of subgroups of PGL$(n)$, 
we refer to, {\it e.g.}, \cite{Sc}.

Let us begin with plane curves relative to SA$(2)$.
Let ${\bf A}^2$ be the unimodular affine plane
with the determinant function $|x\ y|=x^1y^2-x^2y^1$ for two vectors
$x=(x^1,x^2)$ and $y=(y^1,y^2)$,  and let
$x(t)$ be a curve with parameter $t$ into ${\bf A}^2$,
which is nondegenerate in the sense that $|x'\ x''|\neq 0$. When
$|x'\ x''|=1$, furthermore, the parameter $t$ is called 
an equiaffine length parameter.
In this case, 
it holds that $|x'\ x'''|=0$, namely $x'''$ is linearly dependent on
$x'$ and we can write this dependence as
\[x'''=- k_ax',\]
where $ k_a$ is a scalar-valued function 
called the equiaffine curvature. Conversely, given a differential
equation of this form, the map defined by two linearly independent
non-constant solutions defines a curve whose equiaffine curvature is $k_a$.

With reference to this presentation of equiaffine notions, 
we first define 
the general-affine length parameter and the general-affine curvature
relative to the full affine group GA$(2)$ for a nondegenerate curve in
Section \ref{sec:gacurve}. 
In contrast to  the differential equation above, we have
\[x'''=- \frac{3}{2} kx'' - \left(\epsilon + \frac{1}{2}k'+\frac{1}{2} k^2\right)x',
\]
 where $k$ is the general-affine curvature and $\epsilon=\pm 1$ denotes
 additional information of the curve; see Section $\ref{subsectionga}$. 
 Conversely, for a given function $k$ and $\epsilon = \pm1$, there exists
 a nondegenerate curve $x$ for which  
 $x$ satisfies the above equation and the curvature of $x$ is $k$, uniquely
 up to a general-affine motion (Theorem \ref{planenatural}).
We then give remarks on the total curvature of closed curves
and on the sextactic points (Corollary \ref{totalcurv}), and
we study how to compute the curvature. In particular, we give
an expression of the curvature of graph immersions and
classify plane curves with constant general-affine curvature  
(Proposition \ref{gaconst}), and discuss some relations 
with equiaffine treatment and projective treatment of plane curves.

We next consider an extremal problem of the general-affine length
functional and derive a variational formula, 
a nonlinear ordinary differential equation
 that characterizes an extremal plane curve, as
\[ 
k''' + \frac{3}{2} k k'' + \frac{1}{2} {k'}^2 + \frac{1}{2} k^2 k' 
+  \epsilon k' =0
\] 
(Proposition \ref{planeextremalprop}).
Here we give remarks on the preceding studies: the formulation used
to define general-affine curvature of plane curves in this paper is
very similar to that in \cite{Mihai1}. For example, the ordinary
differential equation above for $x$ and Theorem \ref{planenatural}
were already given in \cite{Mihai1}. The formula of the general-affine
plane curvature was given also in \cite{Sc, OST} in a different context.
The variational formula above for $k$ was first given 
in \cite[(33)]{Mihai2}, though some modifications are necessary. 
The same formula for $\epsilon=1$
was then rediscovered by S. Verpoort \cite[p.432]{Ve} in the
equiaffine setting. Furthermore, the author of \cite{Ve}
gave a relation of solutions of this nonlinear
equation with the coordinate functions of the immersion; we 
reprove this relation in Corollary \ref{planeextremal}. 

We then remark that 
the differential equation
above is very similar to some of the nonlinear differential equations
of Chazy type. Furthermore we derive the variational formula
for a certain generalized curvature functional 
 (Theorem \ref{genvarform}).

When the curve is given as a graph immersion, the general-affine
curvature is written as a nonlinear form of a certain intermediate function.
It is interesting to obtain a graph immersion from a given
curvature function, which is treated in Section \ref{sec:graphimmersion}:
We see that its integration can be reduced to solving the Abel 
 equations of the first kind and the second kind (Theorem \ref{abeleq}).

\bigskip

The second aim of this paper is to study general-affine invariants of
space curves.
The procedure is similar to that for the
plane curves. We give the definition
of general-affine space curvatures of two kinds $k$ and $M$,  and 
the ordinary differential equation of rank four
\[
x''''=-3kx'''-\left(2k'+\frac{11}{4}k^2+\epsilon\right)x'' 
  - \left(M+\frac{1}{2}\epsilon k + \frac{1}{2}k''
          +\frac{7}{4}kk' + \frac{3}{4}k^3\right)x'
\]
(Lemma \ref{gaspclemma}), which defines the immersion. Then
we show how to obtain curvatures,
and discuss relations with the equiaffine and projective treatment of
space curves. In particular, we give a list of nondegenerate
space curves of constant general-affine curvature 
and a new proof of the theorem that a nondegenerate curve in the affine
$3$-space is extremal relative to the equiaffine length functional
if and only if the two kinds of equiaffine curvature vanish,
due to \cite{Bl} (Theorem \ref{equi3extremal}).

Finally, we solve
the extremal problem of the general-affine length functional:
a {nondegenerate} space curve without affine inflection point
is general-affine extremal
if and only if the pair of ordinary differential equations 
\begin{align*}
&k''' + \frac{3}{2} k k'' + \frac{1}{2}k^2 k' +\frac{1}{2}
 {k'}^2 
 - \frac{1}{5} \epsilon k' + \frac{6}{5} M'=0,  \\
&k'' + \frac{2}{3} k k' + \frac{5 }{6}\epsilon  k' M 
 -\frac{3 }{2} \epsilon k  M'- \epsilon {M''}=0
\end{align*}
are satisfied (Theorem \ref{spaceextremal}).
Then, we discuss a similarity amongst the nonlinear 
differential equations
for the curvature functions, one for plane curves, and the other
for space curves belonging to a linear complex, 
 {\it i.e.}, $M=\epsilon k$ (Corollary 
\ref{lincomplex}, \ref{planetospace}).

In Appendix, we discuss the projective treatment of plane
curves and space curves and present the variational formula
of the projective length functional by use of the method
in this paper. Theorem \ref{cartanplane} reproduces the
variational formula for projective plane curves due to E.~Cartan,
and Theorem \ref{projspcextremal} gives the variational formula
for projective space curves, which is essentially due to \cite{Ki}.
Furthermore, we treat 
nondegenerate projective homogeneous space curves, called ``W-Kurve''.
The list of such curves may be found elsewhere, but 
nonetheless, we give here a list in Appendix C for later reference.

In this paper, we use the classical moving frame method; refer to
\cite{Ca2, ST}. 
For the equiaffine differential geometry and its terminologies, 
we refer to the books \cite{NS, Bl, Sc}, and, 
for the projective treatment of curves, 
to E. J. Wilczynski \cite{Wi}, E. P. Lane \cite{La}.

\section{General-affine curvature of plane curves} \label{sec:gacurve}

Let $x:M\longrightarrow {\bf A}^{2}$ be a curve
into the $2$-dimensional affine space, where $M$ is a $1$-dimensional
parameter space.
Let $e=\{e_1,e_{2}\}$ be a frame along $x$; at each 
point of $x(M)$ it is a set of independent vectors of ${\bf A}^{2}$
that depends smoothly on the parameter.
The vector-valued $1$-form $dx$ is written as
\begin{equation} \label{dx}
dx = \om^1 e_1 + \om^2 e_2,
\end{equation}
and the dependence of $e_i$ on the parameter is described by the equation
\begin{equation} \label{dei}
de_i = \sum_j \om_i^j e_j,
\end{equation}
where $\om^j$ and $\om_i^j$ are $1$-forms, and the matrix of $1$-forms
\[
\Omega = \left(\begin{array}{cc}
\om^1 &  \om^2 \\
\om_1^1 & \om_1^2 \\
\om_2^1 & \om_2^2
\end{array}
\right)
\]
is called the coframe. 

\subsection{Choice of frames for plane curves and general-affine curvature}\label{frames}
We now reduce the choice of frames in order to define certain invariants.
First, we assume $\om^2=0$, which means that $e_1$ is tangent to the curve, and 
we set $\om^1=\om$ for simplicity. The vector $e_2$ is arbitrary at present,
as long as it is independent of $e_1$. 
Let $\tilde{e}=\{\tilde{e}_1, \tilde{e}_2\}$ be another choice of
such a frame. Then, it is written as
\[\tilde{e}_1 = \lam e_1,\qquad \tilde{e}_2 = \mu e_1 + \nu e_2,\]
where $\lam\nu \neq 0$.
The coframe is written as $\tilde{\om}$ and $\tilde{\om}_i^j$, 
which satisfy
\[ dx=\tilde{\om}\tilde{e}_1, \qquad d\tilde{e}_i = \sum_j \tilde{\om}_i^j\tilde{e}_j.
\]
Then we certainly have
\begin{equation}\label{tildeom}
\tilde{\om} = \lam^{-1}\om.
\end{equation}
Since $d\tilde{e}_1$ is represented in two way, one being
\[ d\tilde{e}_1 = (d\lam) e_1+ \lam (\om_1^1e_1+\om_1^2e_2)\]
and the other being
\[ d\tilde{e}_1 = \tilde{\om}_1^1(\lam e_1)+\tilde{\om}_1^2(\mu e_1+\nu e_2),\]
by comparing the coefficients of $e_1$ and $e_2$ in these expressions, 
we get
\begin{equation}\label{chg1}
\begin{array}{rcl}
& \lam \om_1^2 = \nu \tilde{\om}_1^2,& \\
& d\lam  + \lam \om_1^1 = \lam \tilde{\om}_1^1 + \mu \tilde{\om}_1^2. &
\end{array}
\end{equation}
Similarly, by considering $d\tilde{e}_2$, we have
\begin{equation}\label{chg2}
\begin{array}{rcl}
& \mu \om_1^2+d\nu +\nu \om_2^2 = \nu \tilde{\om}_2^2,& \\
& d\mu  + \mu \om_1^1 +\nu \om_2^1 = \lam \tilde{\om}_2^1 
+ \mu\tilde{\om}_2^2. &
\end{array}
\end{equation}
Since the immersion is $1$-dimensional, we can set
$\om_1^2=h\om$ and $\tilde{\om}_1^2=\tilde{h}\tilde{\om}$, and 
then the first identity of $(\ref{chg1})$ implies that
\[ \tilde{h}=\nu^{-1} \lam^2 h.\]
Hence the property that $h$ is nonvanishing is independent of the frame
and we assume in the following that it is nonvanishing. Such a curve
is said to be \textsl{nondegenerate}.
Geometrically, this property means that the curve is locally strictly convex
at each point. By the identity above, provided that $h$ is nonzero, we can
choose a frame $\tilde{e}$ so that $\tilde{h}=1$ and
we treat such frames with $h=1$ in the following.
Then $\nu=\lam^2$ immediately follows. Next, we see from $(\ref{chg1})$
and $(\ref{chg2})$ that
\begin{eqnarray*}
& \tilde{\om}_1^1 = \om_1^1 + \lam^{-1}d\lam - \mu\lam^{-2}\om, & \\
& \tilde{\om}_2^2 = \om_2^2 + \nu^{-1}d\nu + \mu\nu^{-1}\om. &
\end{eqnarray*}
Hence,
\[ 2\tilde{\om}_1^1-\tilde{\om}_2^2 = 2\om_1^1-\om_2^2 -3\mu\lam^{-2}\om,\]
which means that we can choose $\mu$ so that
$2\tilde{\om}_1^1-\tilde{\om}_2^2 = 0$, and, by considering only such frames
in the following, we must have $\mu=0$.
Thus, we have determined the frame $e$ up to a change of the form
\[ \tilde{e}_1 = \lam e_1, \qquad \tilde{e}_2 = \lam^2 e_2.\]
We call the direction determined by $e_2$ 
the \textsl{general-affine normal direction}. Furthermore, we have
\[\tilde{\om}_1^1 = d\log \lam + \om_1^1, 
\qquad \tilde{\om}_2^1=\lam\om_2^1.\]
From the first identity, we can choose $\lam$ so that $\tilde{\om}_1^1=0$.
Hence, we consider the frame with $\om_1^1=0$ and $\lam$ is assumed to
be constant. From  the second identity, by setting
\[\om_2^1=-\ell \om,\]
and, similarly, $\tilde{\om}_2^1= - \tilde{\ell}\tilde{\om}$, we get
\begin{equation} \label{elltrans}
\tilde{\ell} = \lam^2\ell.
\end{equation}
We call this scalar function $\ell$ the \textsl{equiaffine curvature}, see Section 
\ref{equigen}, or \textsl{affine mean curvature}, 
in analogy with equiaffine theory of hypersurfaces, though it still
depends on the frame chosen. A point where $\ell=0$ is called
an \textsl{affine inflection point}.
For its geometrical meaning, we refer to Section $\ref{subsc:graph}$
and \cite{IS}.
Thus, we have seen that, given a nondegenerate curve $x$, there exists a frame
$e$ with coframe of the form 
\begin{equation}\label{eq:coframeforequi}
\left(
\begin{array}{cc}
\om & 0 \\
0   & \om \\
-\ell \om & 0
\end{array}\right)
\end{equation}
and that such frames are related by $\tilde{e}_1=\lam e_1$ 
and $\tilde{e}_2=\lam^2 e_2$ for a nonzero constant $\lam$.
In the following, given a curve $x=x(t)$ with parameter $t$,
we assume that the vector $e_1$ is a positive multiple of the
tangent vector $dx/dt$. Then, the choice of $\lam$ is 
limited to be positive and the form $\om$ is
a positive multiple of $dt$.

We now assume $\ell\neq 0$ and 
let $\epsilon$ denote the sign of $\ell$: 
\[\epsilon={\rm sign}(\ell).\]
 It is a locally defined invariant of the curve called the \textsl{sign} 
 of the curve.
Then we define a form 
\begin{equation}
\om_s = \sqrt{\epsilon \ell}\om,
\end{equation}
which is uniquely defined, independent of the frame, in view of 
$(\ref{tildeom})$ and $(\ref{elltrans})$.
We call this form the \textsl{general-affine length element} and
call the parameter $s$ such that $ds=\om_s$
the \textsl{general-affine length parameter}, 
determined up to an additional constant.
\begin{definition} \label{planecurv}
We call the scalar function $k$ defined as
\[ \frac{d\ell}{\ell} = k\om_s,\]
the {\rm general-affine curvature}.
In other words,
\begin{equation}
k={d\log\ell\over ds}.
\end{equation}
\end{definition}
We define a new frame $\{E_1, E_2\}$ by setting
\begin{equation}\label{eq:E1E2}
 E_1 = \frac{1}{\sqrt{\epsilon\ell}}e_1,\qquad E_2=\frac{1}{\epsilon \ell}e_2.
\end{equation}
Then, 
\[dx = \om_s E_1.\]

For another frame $\{\tilde{e}_1,\tilde{e}_2\}$
where $\tilde{e}_1=\lambda e_1$ and $\tilde{e}_2 = \lambda^2 e_2$, 
 we similarly define
 $\tilde{E}_1$ and $\tilde{E}_2$. Then we can see that
\[
 \tilde{E}_1={1\over \sqrt{\epsilon\tilde{\ell}}}\tilde{e_1}
={1\over \sqrt{\epsilon \lambda^2\ell}}\lambda e_1 = E_1\]
and
\[\tilde{E}_2={1\over \epsilon\tilde{\ell}}\tilde{e_2}
  ={1\over \epsilon \lambda^2\ell}\lambda^2 e_2 = E_2. \quad \;\;\]
Thus, we have proved the following.

\begin{proposition} 
Assume $\ell\neq 0$. Then, the frame $\{E_1, E_2\}$ 
is uniquely defined from the immersion and it satisfies a Pfaffian equation
\begin{equation} \label{pfaff}
 d \left(\begin{array}{c} x \\ E_1 \\ E_2\end{array}\right)
= \Omega
\left(\begin{array}{c} E_1 \\ E_2 \end{array}\right);
\qquad \Omega =
\left(
\begin{array}{cc}
\om_s & 0 \\
-{1\over 2}k\om_s   & \om_s \\
-\epsilon \om_s & -k\om_s
\end{array}\right),
\end{equation}
where $\om_s$ is the general-affine length form, 
$k$ is the general-affine curvature and $\epsilon$ is ${\rm sign}(\ell)$. 
\end{proposition}

By use of this choice of frame, we have the following lemma.
\begin{lemma}\label{lem:kode} The immersion $x$ satisfies the ordinary differential
equation
\begin{equation} \label{kode}
x'''+{3\over 2}kx'' + \left(\epsilon + {1\over 2}k'+{1\over 2}k^2\right)x'=0,
\end{equation}
relative to a general-affine length parameter.
\end{lemma}

\noindent Proof. The equation $(\ref{pfaff})$ shows that
$x'=E_1$, $E_1'=-{1\over 2}kE_1 + E_2$, and $E_2'=-\epsilon E_1-kE_2$,
where the derivation $\{{}'\}$ is taken relative to the length parameter. 
Then, combining these derivations, we easily obtain
the differential equation above.

\begin{remark} \label{paramsign}
{\rm The definition of the curvature depends on the orientation
of the parameter $t$.
If we let the parameter be $u=-t$ and denote by an overhead 
 dot the derivation
relative to $u$, then we have
\[ 
\threedot{x}-{3\over 2}k \ddot{x} + \left(\epsilon - {1\over 2}\dot{k}
+{1\over 2}k^2\right)\dot{x}=0.
\]
 Namely, the curvature changes sign and its absolute value is 
a true invariant independent of the orientation of the parameter.}
\end{remark}

With this remark in mind, we have the following theorem.

\begin{theorem}[\cite{Mihai1}] \label{planenatural}
Given a function $ k(t)$ of a parameter $t$ and $\epsilon=\pm 1$, 
there exists a nondegenerate curve $x(t)$ for which $t$ is
an length-parameter, $k$ the curvature function and $\epsilon$ the sign
of $\ell$,  uniquely up to a general-affine transformation.
\end{theorem}

\noindent Proof. Given $k$ and $\epsilon$, 
we solve the ordinary differential equation in \eqref{kode}
to get the vector $x'(t)$, 
which is determined up to a  general linear transformation.
Then, we get $x(t)$ up to an additional translation by 
a constant vector; that is, the curve $x(t)$ is determined 
up to a transformation in GA$(2)$.  

Theorem \ref{planenatural} and the ordinary
differential equation \eqref{kode}
were first given by T. Mih$\breve{\rm a}$ilescu
in \cite{Mihai1}, to the authors' knowledge; refer
also to \cite{Mihai2} and \cite{CG}.

\begin{example} \label{gacircle}
Ellipse and Hyperbola.
{\rm
Let $x$ denote an ellipse $(a\cos\theta,b\sin\theta)$ or
a hyperbola $(a\cosh\theta, b\sinh\theta)$.  Then, $x'''=-\epsilon x'$,
where $\epsilon=1$ for the ellipse and $\epsilon=-1$ for
 the hyperbola. 
It is easy to see that $\theta$ is a general 
 affine length, see \eqref{dstwo2}.
 Hence, $k=0$. }
\end{example}
 According to this example, we may call a nondegenerate curve is
 of \textsl{elliptic {\rm (resp.} hyperbolic{\rm )} type} if $\epsilon=1$
 (resp. $\epsilon=-1$).

 We call the vector $E_2$, uniquely defined when $\ell\neq 0$, 
the \textsl{general-affine normal} and the map 
$t \longmapsto E_2$ the general-affine Gauss map. 
Then, by an analogy with affine spheres in equiaffine differential geometry,
it is natural to call a curve such that the map $E_2$ passes
through one fixed point a \textsl{general-affine circle}. 
For the ellipse or the hyperbola, $E_2= x'' $ and it holds that
\[x+\epsilon E_2=0,\]
and thus it is a general-affine circle. Conversely,
for a curve $x$ to be such a circle, there exists a scalar function $r(t)$
and a fixed vector $v$ such that
\[ x+ rE_2=v.\]
However, this implies that $dx+  dr E_2  + rdE_2=0$, which induces,
by the identity $(\ref{pfaff})$, $(1-\epsilon r)\om E_1 + (dr-kr\om)E_2=0$.
Hence, $r=\epsilon$ is constant and $k=0$.
Then, by integrating the differential equation $(\ref{kode})$ when $k=0$,
we see that
any general-affine circle is general-affinely congruent to (a part of)
an ellipse or a hyperbola. We also refer to Example \ref{ellzero}.

\subsection{Total curvature and sextactic points}

The formula $(\ref{pfaff})$ implies the identity
\begin{equation} \label{traceOm}
d\log\left(\left|\det\left(\begin{array}{cc} E_1\\E_2\end{array}\right)
\right|\right)
 = -{3\over 2}k\om_s,
\end{equation}
where $\det$ is taken relative to a (any) unimodular
structure of the space ${\bf R}^2$.  This formula shows the
following corollary immediately.

\begin{corollary} \label{totalcurv}
Assume that the curve $C$ is nondegenerate and closed, and 
has no affine inflection point.
Then, the total curvature
$\int_C  k\om_s$ vanishes. In particular, such a curve has at least
two general-affine flat points.
\end{corollary}

As we will see in Section \ref{gatoproj}, any general-affine
flat point, where $k=0$ by definition, is nothing but a sextactic point.
We know a classical
theorem due to Mukhopadhayaya, also due to G. Herglotz and J. Rado, 
we refer, {\it e.g.},  to  \cite{ST,TU},
that the number of sextactic points of a strictly convex simply
closed smooth curve is at least six. In other words, on such a
curve there are at least six general-affine flat points.

Furthermore, as an analogue of the Euclidean plane curve, it is
natural to introduce a notion of a \textsl{general-affine vertex} where
$k$ is extremal. The corollary above says that any nondegenerate closed curve
without affine inflection point has at least two general-affine
vertices. In fact, 
the example \ref{rose2} shows that there exists a plane curve which 
has two general-affine vertices.


\subsection{Computation of general-affine curvature of plane curves}
\label{subsectionga}

In this subsection, we will see how to obtain the curvature of a
curve given relative to a parameter not necessarily a length parameter.

Let $t\longrightarrow x=x(t)\in {\bf A}^2$ be a nondegenerate curve so that
the vectors $x'$ and $x''$ are linearly independent. Then the derivative
$x'''$ is written as a linear combination of $x'$ and $x''$: there are
scalar functions $a=a(t)$ and $b=b(t)$ such that
\begin{equation} \label{xode}
x''' = a x'' + b x'.
\end{equation}
Since $dx=x'\, dt$, the frame vector $e_1$ is a scalar multiple of $x'$:
\begin{equation}\label{eq:frame}
 dx = \om \, e_1;\qquad e_1=\lambda x',\quad \om = \lambda^{-1}dt.
\end{equation} 
Then, the derivation 
\[ de_1=\lambda ( \lambda x'' + \lambda' x')\,  \om\]
implies  that the second frame vector is
\[ e_2=\lambda ( \lambda x'' + \lambda' x').\]
In order for the frame $\{e_1, e_2\}$ to be chosen as in 
Section \ref{frames}, the vector $de_2$ must be a multiple of $e_1$. Since
\begin{eqnarray*}
de_2 &=& \lambda(\lambda^2x'''+3\lambda\lambda'x''+
   (\lambda\lambda''+{\lambda'}^2)x')\, \om \\
  &=& \{\lambda^2(\lambda a + 3\lambda')x'' +
        \lambda(\lambda^2b+\lambda\lambda''+{\lambda'}^2)x'\}\, \om,
\end{eqnarray*}
we have
\begin{equation}\label{eq:lambda}
\lambda a + 3\lambda'=0, \quad {\it i.e.}\quad
   \lambda =  e^{-{1\over 3}\int a(t) dt}
\end{equation}
 up to a positive constant multiple, and by definition, 
\begin{equation}\label{eq:ell}
\ell = - (\lambda^2b+\lambda\lambda''+{\lambda'}^2).
\end{equation}
We now assume that $\ell\neq 0$ and 
recall that $\epsilon = {\rm sign} (\ell)$.
Then, we have
\[
ds^2 = -\epsilon\left(
 b+{\lambda\lambda''+\lambda'^2\over \lambda^2}\right) dt^2.
\]
In terms of $a$ and $b$, 
\begin{equation}\label{dstwo2}
ds^2 =-\epsilon \left(b+{2\over 9}a^2 - {1\over 3}a' \right) dt^2.
\end{equation}
Hence, a length parameter $s$ which is a function of $t$ is
obtained by solving the equation
 \[ \left({ds \over dt}\right)^2 
   =-\epsilon \left(b+{2\over 9}a^2 - {1\over 3}a'\right).\]
If in particular $t$ itself is a length parameter, then we must have
\begin{equation}\label{length}
\ell=\lambda^2\epsilon, \qquad 
 -b - {2\over 9}a^2 + {1\over 3}a'=\epsilon.
\end{equation}

Assume now that the curve $x(t)$ is given relative to a length parameter $t$.
Then the differential equation $(\ref{xode})$ is written as
\begin{equation} \label{gaode}
x''' = a x'' + \left({1\over 3}a'-{2\over 9}a^2 -\epsilon \right) x',
\end{equation}
and the curvature is
\begin{equation} \label{odecurv}
k  = {d\log (\epsilon\lambda^2 ) \over dt} = -{2\over 3}a,
\end{equation}
which agrees with the expression of
the coefficient in the equation $(\ref{kode})$.
\medskip

For another parameter 
$\sigma=\sigma(t)$, 
we write
\[ y(\sigma) = x(t).\]
Then, a calculation shows that
\[ \threedot{y}(\sigma) = A(\sigma)\ddot{y}(\sigma)+B(\sigma)\dot{y}(\sigma),\]
where
\begin{eqnarray}
&A(\sigma) = \displaystyle
\left(a-{3\sigma''\over \sigma'}\right){1\over \sigma'},& \label{newa}\\
&B(\sigma) = \displaystyle  \left(b + a{\sigma''\over \sigma'} 
          - {\sigma'''\over \sigma'}\right){1\over \sigma'^2}.&
\label{newb}
\end{eqnarray}
We note that there holds a covariance relation:
\[ B+{2\over 9}A^2-{1\over 3}\dot{A}
  =\left(b+{2\over 9}a^2-{1\over 3}a'\right){1\over \sigma'^2}.\]
Thus, we have the following procedure to obtain curvature:
\bigskip

\noindent {\sf Procedure for computing curvature
\begin{enumerate}
\item Given a curve $x(t)$, derive the differential equation $(\ref{xode})$.

\item {C}ompute $L=\left(-b-{2\over 9}a^2 + {1\over 3}a' \right)$ and
define $\epsilon$ by $\epsilon={\rm sign}(L)$.

\item {C}ompute the length parameter $\sigma$
by solving $d\sigma = \sqrt{\epsilon L}dt$.

\item {C}ompute $A$ by $(\ref{newa})$; then, $-{2\over 3}A$ is the curvature.
\end{enumerate}
}

\begin{example} \label{logsp}
{\rm A logarithmic spiral is the curve
$x(t)=e^{\gamma t}(\cos \alpha t, \sin \alpha t)$. It is easy to see
that 
$x''' = 2\gamma x'' - (\gamma^2+\alpha^2)x'$. Hence, $a=2\gamma$,
$b=-\gamma^2 -\alpha^2$, and $L=\gamma^2/9+\alpha^2$;
hence, $\epsilon=1$.  The length parameter is $s=\sqrt{\gamma^2/9+\alpha^2}t$
and, by rewriting the equation with this $s$, the coefficient $a$ is multiplied
by $1/\sqrt{\gamma^2/9+\alpha^2}$. Hence, the curvature is the
constant $-4\gamma/\sqrt{\gamma^2+9\alpha^2}$. We require here 
 that $\gamma\neq 0$, since the curve when $\gamma=0$ is a circle,
 which we will consider in Example $\ref{ellzero}$, 
 and also $\alpha\neq 0$, because  the curve is nondegenerate. 
 Then, it is easy to see that
the possible values of $k$ range in $0<|k|<4$. 
}
\end{example}

\begin{example} \label{catenary}
{\rm The catenary curve is defined as $x(t)=(t,\cosh(t))$. The equation is
\[x'''=\tanh(t)x''.\]
Then $L=-(2\cosh(t)^2-5)/(9\cosh(t)^2)$, which 
vanishes at $2\cosh(t)^2-5=0$ ($t=\pm 1.031\dots$) , 
where the curvature is undefined.
For the value $t$ where $L<0$, on the dotted curved of the left-hand 
side
picture in Figure 1, $\epsilon=1$ and, elsewhere $\epsilon=-1$.
As the picture shows, it is difficult to ``see'' where $L$ vanishes
and how $\epsilon$ changes the value.}
\end{example}

\begin{example} \label{rose2}
{\rm The curve $x(t)=(\cos(nt)\cos(t), \cos(nt)\sin(t))$ is called
a rose curve. Now choose $n=1/3$ with the range $t\in [0,3\pi]$; 
shown in Figure 1 (right).
It satisfies the equation
\[ x''' = -{8\sin(t/3)T\over 3(1+4T^2)}x''
        -{4(7+8T^2)\over 9(1+4T^2)}x',\]
where $T=\cos(t/3)$. 
A computation shows that $L>0$, $\epsilon=1$ and 
the length parameter $\sigma$ is defined by
$d\sigma = 
{2 \sqrt{256T^2+320T^4+69}\over 9(1+4T^2)}dt$,
and the curvature $k(t)$ is defined for all value $t$. It 
vanishes at $t=0, 3\pi/2$.
Moreover, it is easy to see that 
 the number of general-affine vertices is two. 
Thus the rose curve with $n=1/3$ attains
 the minimum number of general-affine vertices 
amongst general-affine plane closed curves.
}
\end{example}

\begin{figure}
\centering
\begin{tabular}{cc}
\includegraphics[width=6cm]{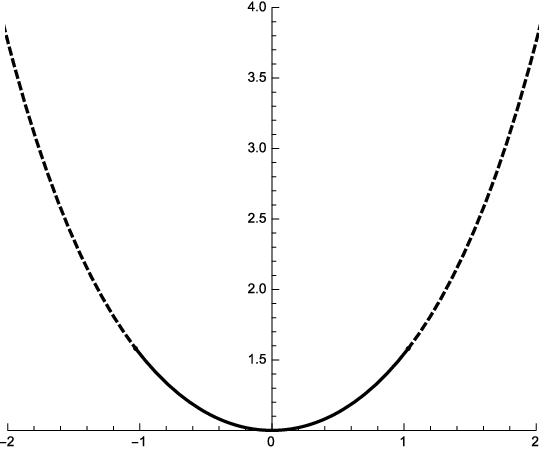} \qquad & \qquad 
\includegraphics[width=6cm]{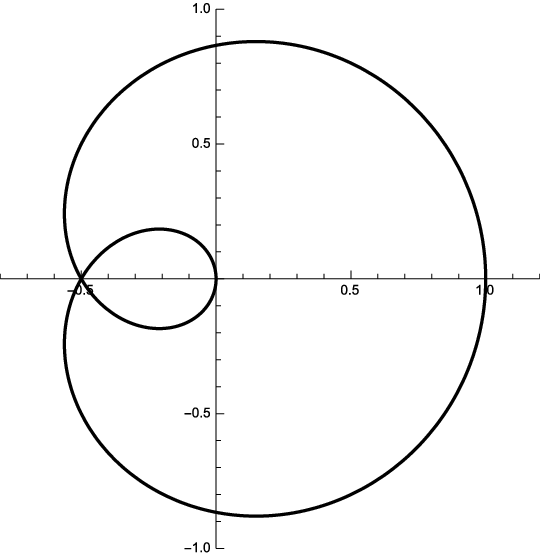} \\
\multicolumn{2}{c}{Figure 1. {\rm Catenary (left) and rose curve (right)}}
\end{tabular}
\end{figure}

\comment{
\begin{example} \label{rose3}
{\rm The curve $x(t)=(\sin(3t)\cos(t), \sin(3t)\sin(t))$ is one of
rose curves. It satisfies the equation
\[x''' = -{24 \cos(3t)\sin(3t)\over 4\cos(3t)^2+5}x''
-{4(35-8\cos(3t)^2)\over 4\cos(3t)^2+5}x'.\]
A computation shows that $\epsilon=1$ and 
the length parameter $\sigma$ is defined by
$d\sigma = {2\sqrt{205-16\cos(3t)^2}\over 4\cos(3t)^2+5}dt$ and, that
the curvature $k(t)$ is defined for all value $t$, and has a
symmetry $k(\pi/3-t)=-k(t)$. 
}
\end{example}
}

\subsection{General-affine curvature of a graph immersion}\label{subsc:graph}
 Let us consider the nondegenerate curve given by a graph immersion
 $x(t)=(t, f(t))$.
 We will find the formula of the curvature given by the function $f$ and 
 show fundamental examples of graph immersions.

Note that since $x$ is nondegenerate, $x''=(0, f'') \neq 0$, we 
 can assume $f''>0$  in the following.
 Since $x'''=(0,f''')$,
the coefficients of the differential equation $(\ref{xode})$ are
\[ a={f'''\over f''}, \qquad b=0.\]
Hence,
\[\lambda = e^{-\frac{1}{3} \int a(t) dt} =  (f'')^{-1/3}\]
up to a constant multiple. 
With this $\lambda$, we have
\[ \ell = -(\lambda\lambda'' + (\lambda')^2)\]
and the length element is
\[ ds^2 = -\epsilon \lambda^{-2}(\lambda\lambda'' + (\lambda')^2) dt^2.\]
If we set $\mu=\lambda^2= (f'')^{-2/3}$, then 
\begin{equation} \label{muell}
\ell=-\frac{\mu''}{2},\qquad ds^2 = -{\epsilon \mu''\over 2\mu}dt^2.
\end{equation}
Hence, we have the formula
\begin{equation} \label{curv}
 k = {d\log \ell\over ds} 
= \sqrt{-2\epsilon\mu\over \mu''}{\mu'''\over \mu''}.
\end{equation}

\begin{lemma} The quantities $ds^2$ and $k^2$ are expressed by
use of the function $f$ as follows$:$
\begin{equation} \label{olver}
\begin{array}{c}
\displaystyle 
ds^2 = {\left|3f''f'''' - 5(f''')^2 \right|\over 9(f'')^2} dt^2,\\
\noalign{\medskip}
\displaystyle
k^2 = { \left|9(f'')^2f''''' - 45f''f'''f'''' + 40(f''')^3\right|^2
   \over \left|3f''f'''' - 5(f''')^2\right|^3}.
\end{array}
\end{equation}
\end{lemma}
This is seen by expressing $\mu''$ and $\mu'''$ explicitly in terms of 
derivations of $f$. As we will see in Section $\ref{equigen}$,
$-\mu''(= 2 \ell)$ equals the equiaffine curvature up to a multiplicative
constant. The factor in the right-hand side of the first expression
is known, see \cite[p. 14]{Bl}.
The second expression of $k^2$ was already presented in \cite[p. 54]{Sc},
which was proved from a different point of view. We refer also to
\cite[p. 343]{OST} and \cite{CQ}. The differential
polynomials in the numerator and
 denominator in the right-hand side of $k^2$ are known
classically; Berzolari \cite{Be} stated that those go back to
G. Monge in 1810: Sur les \'Equations diff\'erentielles des
courbes du second degr\'e, Corresp. sur l'\'Ecole imp.
Polytechn., Paris, N$^{\rm o}$ II, 1810, pp. 51--54.

\vskip1pc

We consider a nondegenerate curve around the point $t=0$.
For appropriate affine coordinates, we have an expression
\[ {f} = {1\over 2}t^2 + {p\over 4!}t^4 + {q\over 5!}t^5 + \cdots.\]
Then, we see
\[ \mu = 1- {p\over 3} t^2 - {q\over 9}t^3 + \cdots\]
and
\begin{equation} \label{muexp}
 \mu'' = -{2p\over 3} - {2q\over 3}t + \cdots.
\end{equation}
Hence, the property $\ell=0$ at $t=0$
means that $p=0$ and, hence, that the osculating
parabola touches the curve to at least $5$-th order.
In particular, $\ell\equiv 0$ for any parabola.
Conversely, we have the following example.

\begin{example}
 \label{ellzero}
\textit{Curves with constant $\ell$}.
{\rm We first let $\ell\equiv 0$. Then, $\mu=at+b$ or $\mu =a$,
and it follows that $f=c_1(at+b)^{1/2}+c_2t+c_3$ 
or $f=c_1t^2+c_2t+c_3$ for certain constants $a, b, c_1, c_2, c_3$.
Namely, the curve $x(t)= (t, f(t))$ is ({general-}affinely equivalent to) a parabola.

If $\ell$ is a nonzero constant: $\ell=-c \neq 0$, 
then  $\mu=ct^2+at+b$, and 
$f(t)=c_1(ct^2+at+b)^{1/2}+c_2t + c_3$ for constants $a,b,c,c_1,c_2,c_3$.
This implies that the curve is an ellipse or a hyperbola.
Thus, we have seen that the curve with constant  $\ell$ is a quadric
 and hence the curvature satisfies $k=0$. 

}
\end{example}
 For the ellipse ${f}=-(\alpha^2 - t^2)^{1/2}$ and the hyperbola 
 ${f}=(\alpha^2+ t^2)^{1/2}$,
 we can see that $\mu=\alpha^{-4/3}(\alpha^2\pm t^2)$
 and $\displaystyle ds^2 = {1\over \alpha^2\pm t^2}dt^2$.  
Relative to the reparametrization by use of angular variable
$t=\alpha\sin\theta$ or $\alpha\sinh \theta$,
we have $ds^2=d\theta^2$; namely, $ds$ is independent of $\alpha$, 
the ``size'' of the curve in euclidean sense.

\comment{
 The vector $E_2$ in \eqref{eq:E1E2} is uniquely defined when $\ell\neq 0$, which
 we may call the \textsl{general-affine normal}, and the map $t\rightarrow
 E_2$ is the general-affine Gauss map. 
 Then, it is natural to call a curve such that the map $E_2$ passes
through one fixed point a \textsl{general-affine circle}. 
For the ellipse or the
hyperbola ${f}= -\epsilon(\alpha^2 - \epsilon t^2)^{1/2}$, it holds that
\[x+\epsilon E_2=0\]
 and thus it is a general-affine circle.
 Conversely for a curve $x$ to be such a general-affine circle, 
 there exists a scalar function $r(t)$
and a fixed vector $v$ such that
\[ x+ rE_2=v.\]
However, this implies that $dx+  dr E_2  + rdE_2=0$, which induces,
by the identity $(\ref{pfaff})$, $(1-\epsilon r)\om E_1 + (dr-kr\om)E_2=0$.
Hence, $r=\epsilon$ is constant and $k=0$, which implies that
any general-affine circle is general-affinely congruent to (a part of)
an ellipse or a hyperbola.
}

\begin{example} \label{exp} {\rm
For the curve ${f}=e^t$, we see that $\epsilon = -1$ and 
$k=-\sqrt{2}$.
For the curve $f = t\log t$ ($t>0$), we have $\epsilon=1$ and $k=-4$.}
\end{example}

\begin{example} \label{power}
{\rm For the curve ${f}=t^\alpha$, we see that the curvature is
\[ k(\alpha) = {-2(\alpha +1)\over \sqrt{|(2\alpha-1)(\alpha -2)|}}.\]
Here we assume $\alpha \neq 0,\  \pm 1,\  1/2,\  2$ so that the curve is
neither trivial nor quadratic. Since 
$\mu''=(1/9)(2\alpha-4)(2\alpha-1)t^{-2(\alpha+1)/3}$, we have 
$\epsilon =1$ when
$1/2<\alpha<2$, and $\epsilon=-1$ when $\alpha>2$ or $\alpha<1/2$.
Note that $k(1/\alpha)=k(\alpha)$ when $\alpha>0$, and
$k(1/\alpha)=-k(\alpha)$ when $\alpha<0$. By this symmetry, 
in order to know the possible values of $k=k(\alpha)$, it is sufficient 
to consider the case $-1<\alpha<1$. Then, 
$\epsilon=-1$ and $k \in (-\infty,-\sqrt{2})$ for $\alpha\in (0,1/2)$;
$\epsilon=-1$ and $k \in (-\sqrt{2},0)$ for $\alpha\in (-1,0)$;
$\epsilon=1$ and $k \in (-\infty,-4)$ for $\alpha\in (1/2,1)$.
}
\end{example}

Due to Theorem $\ref{planenatural}$, the curves
with constant curvature are obtained by solving the equation 
$(\ref{gaode})$ for constant $a$:
\[
x''' = a x'' + \left(-{2\over 9}a^2 -\epsilon \right) x'.
\]
The case where $a=0$ is the case $\ell=0$: Example $\ref{ellzero}$.
When $a\neq 0$ we have the following:
\begin{proposition} \label{gaconst}
When the curvature $k$ is a nonzero constant,
the curve is general-affinely congruent to one of the curves
in Examples $\ref{logsp}$, $\ref{exp}$ and $\ref{power}$:
$e^{\gamma t}(\cos t, \sin t)$, $(t, e^t)$, $(t, t\log t)$,
$(t, t^\alpha)\; (\alpha \neq 0,\ \pm 1,\  2,\  1/2)$.
\end{proposition}

\noindent Proof. We set $u=e^{-at/2}x'$. Then $u$ satisfies
the equation $u'' + p u=0$, $p=\epsilon - a^2/36$. According to
$p=0$, $p<0$, $p>0$, a set of independent solutions gives
a map $x$ which is congruent to the curves $(t, t\log t)$,
$(t, t^\alpha)$ and $(t, e^t)$, and $e^{\gamma t}(\cos t, \sin t)$
with parameter renewed appropriately.
More precisely, when $\epsilon=1$, we have the curve 
$(t, t^\alpha)$ ($1/2<\alpha<1$), 
$(t, t\log t)$, and  $e^{\gamma t}(\cos t, \sin t)$ 
according to $k\in (-\infty,-4)$, $k=-4$, and $k \in (-4,0)$,
respectively.
When $\epsilon=-1$, we have $(t, t^\alpha)$ ($0<\alpha<1/2$), 
$(t, e^t)$, and $(t, t^\alpha)$ ($-1<\alpha<0$) according to
$k\in (-\infty,-\sqrt{2})$, $k=-\sqrt{2}$, and $k \in (-\sqrt{2},0)$, 
respectively. 

We remark that any curve with constant curvature is an orbit of
a $1$-parameter subgroup of GA$(2)$, because of the unique and 
existence in Theorem
$\ref{planenatural}$.

 Table {\ref{table:gap}} is the 
 classification of general-affine curves with constant curvature $k$. 
  We let $k\leq 0$, see Remark \ref{paramsign}, and note that the case
 $k = -\infty$ corresponds to the parabola defined in Example \ref{ellzero}. 
\bigskip 

\begin{table}
\centering

\begin{tabular}{lll }
\hline
 curvatures 
& 
 curves with $\epsilon =+1$&  Examples \\
\hline
$ k=0$ & $(t, - (\alpha^2 -t^2)^{1/2})$ & \ref{ellzero}\\ 
$-4 < k <0$ & $e^{\gamma t} ( \cos \alpha t, \sin \alpha t)\;\; (\gamma \neq 0, \alpha \neq 0)$ & \ref{logsp} \\ 
$k=-4$ & $(t, t \log t)\; (t>0)$ &  \ref{exp} \\ 
$-\infty < k<-4$  & $(t, t^{\alpha})$ $(\alpha \in (1/2, 1))$&  \ref{power}\\
 \hline  
\end{tabular}
\medskip

\begin{tabular}{lll}
\hline
 curvatures 
& 
 curves with $\epsilon =-1$ &  Examples \\ \hline
$k=0$ &$(t, (\alpha^2 + t^2)^{1/2}$) & \ref{ellzero} \\ 
$k = - \sqrt{2}$ & $(t, e^t)$&  \ref{exp} \\ 
  $-\infty <k < 0, \; k \neq -\sqrt{2}$ & $(t, t^{\alpha})$  
  $(\alpha \in (0, 1/2)$ or $(\alpha \in (-1, 0))$ &  \ref{power}\\ \hline
\end{tabular}
 \caption{Plane curves with constant general-affine curvature}
\label{table:gap}
\end{table}

\subsection{From equiaffine to general-affine} \label{equigen}
Since the group of equiaffine motions
${\rm SA}(2)={\rm SL}(2,{\bf R})\ltimes {\bf R}^2$
is a subgroup of the general-affine group 
${\rm GA}(2)={\rm GL}(2,{\bf R})\ltimes {\bf R}^2$, any general-affine
invariant is obviously an equiaffine invariant. In this subsection,
we give the expression of the general-affine length parameter and 
the general-affine curvature 
by use of the equiaffine length parameter and the equiaffine curvature.

Let us consider the coframe \eqref{eq:coframeforequi}:
\begin{equation}\label{eq:coframeforequi2}
d e = 
\left(
\begin{array}{cc}
\om & 0 \\
0   & \om \\
-\ell \om & 0
\end{array}\right) e
\end{equation}
 for $e = \{e_1, e_2\}$.
 In the equiaffine treatment, it is enough to consider only the unimodular change of 
 frame, {\it i.e.} $\lam \nu =1$. Because $\nu = \lam^2$, we have $\lam^{3} =1$.
 Thus $\ell$ is an absolute invariant. The scalar $\ell$ is usually denoted by $k_a$
 and called the \textsl{equiaffine curvature} of a plane curve.
 The parameter $t$ for which $\omega = dt$ holds is called the 
 \textsl{equiaffine length
 parameter}. Let $x(t)$ be a curve with equiaffine length parameter $t$. Then, 
 it is easy to see by \eqref{eq:coframeforequi2} that $x$ satisfies
\begin{equation} \label{equiaffeq}
 x''' = -k_a x'.
\end{equation}
On the other hand, when the curve is given by a graph immersion 
$x(t)=(t,f(t))$,
where $t$ is a parameter that is not necessarily equiaffine, 
the equiaffine length element is $(f'')^{1/3} dt$ and
the equiaffine curvature
is $k_a = -{1\over 2}\mu''$ for $\mu = (f'')^{-2/3}$.
As we have seen in $(\ref{muell})$, the equiaffine curvature
is nothing but $\ell$ up to a constant multiple.
We refer to \cite[p. 13--14]{Bl}. 

Now we consider the curve in view of the group GA$(2)$. Then, 
the differential equation $(\ref{equiaffeq})$ above shows 
that the general-affine length element is
\[ 
ds = |k_a|^{1/2} dt.
\]
We rewrite the differential equation using the parameter $s=s(t)$:
we set $y(s)=x(t)$ and let $\{\dot{}\}$ denote the derivation relative to $s$. 
Then, we get the equation:
\[ \threedot{y} = A(s) \ddot{y} + B(s)\dot{y},\]
where $A(s)$ and $B(s)$ are given by use of \eqref{newa} and \eqref{newb}. 
For simplicity,
we set $K=|k_a|$. Since $s'=ds/dt = K^{1/2}$ and
$s''={1\over 2}K' K^{-1/2}$, we see that
$A(s) = -{3\over 2}K'K^{-3/2}$. 
Therefore, the general-affine curvature of $x(t)$ is 
\[ k= K'{K^{-3/2}}.\]
The quantity of this form was already treated in \cite[p. 24]{Bl} 
by dimension considerations to get an invariant relative to
similarity transformation. There a remark was given that
the curves with constant ${K'K^{-3/2}}$ consist of parts of curves
called ``$W$-Kurve'', discussed by Klein and Lie \cite{KL}.
\bigskip

\subsection{From general-affine to projective}  \label{gatoproj}

It was G. H. Halphen \cite{Ha} who began
a systematic study of projective curves
in view of ordinary differential equations. Later, E. J. Wilczynski 
gave a  classical treatment of projective curves in the book \cite{Wi}.
Also, the books by E. P. Lane \cite{La} are standard references
for this subject. 
In this subsection, we recall the definition of the
projective length element and the projective curvature, and
give the expressions of such invariants in terms of
general-affine invariants.

A nondegenerate curve in ${\bf P}^2$ with parameter $t$ is described
by an ordinary differential equation of the form
\begin{equation} \label{projone}
  y''' + P_2 y' + P_3 y=0,
\end{equation}
such that a set of three independent solutions, say, $x^1$, $x^2$, $x^3$
defines a map $t\rightarrow [x^1,x^2,x^3]\in {\bf P}^2$, where
$[\quad ]$ denotes homogeneous coordinates. 
For this equation, the form
\begin{equation} \label{projlen}
  P^{1/3}dt,\quad {\rm where} \quad P=P_3 -{1\over 2}P_2',
\end{equation}
is called the projective length element.
Furthermore, when $t$ itself is a projective length parameter, the
equation can be transformed by a certain change of variables from $y$
to $z=\lambda y$ into the equation of the form
\begin{equation} \label{halphen}
  z''' + 2k_p z' + (1+k_p')z=0,
\end{equation}
which is called the Halphen canonical form. 
Then, the coefficient $k_p$ is called the \textsl{projective curvature}
and is given by the formula
\begin{equation} \label{p2curv}
k_p=P^{-2/3}\left({1\over 2}P_2 - {1\over 3}{P''\over P} 
+ {7\over 18}\left({P'\over P}\right)^2\right);
\end{equation}
we refer to \cite[p. 71]{Ca2}.
In particular, when $k_p$ is constant, the curve is called an
anharmonic curve and is obtained by integrating the differential
equation $z''' + 2k_p z' + z =0$; we refer to \cite[p. 86--91]{Wi}.
The list of anharmonic curves is 
the same as the list of plane curves with constant
general-affine curvature in Section \ref{subsc:graph}, up to
projective equivalence; we note that the sign
$\epsilon$ plays no role in the projective classification.

In the general-affine setting we had a differential equation
$(\ref{gaode})$, which can be transformed into the equation 
of the form $(\ref{projone})$ by changing $x$ into 
$y=e^{-{1\over 3}\int a dt}x$. The result is
\[
y''' = \left({1\over 9}a^2 -{2\over 3}a' - \epsilon \right) y'
+\left({1\over 9}aa' - {1\over 3}a'' - {1\over 3}a\epsilon\right)y.\]
Hence, we can see
\[ P= -{1\over 3}a\epsilon;\]
this implies that the projective length element is
$a^{1/3}dt$ up to a scalar multiple, while $-{2\over 3}a$ is the
general-affine curvature.
In particular, the point where the general-affine curvature vanishes
is the point where the invariant $P$ vanishes, which is classically
 called a \textsl{sextactic point}.

We remark that, in \cite{SS}, S. Sasaki showed how to obtain
the projective length parameter and the projective curvature
directly from the equiaffine curvature.

\section{General-affine extremal plane curves and the associated differential 
  equation}

 In Section $\ref{frames}$, we have defined the general-affine length element 
 $\om_s =\sqrt{\epsilon \ell} \omega$ with $\epsilon = {\rm sign} (\ell)$.
 It defines the length functional for general-affine curves, 
and in this section 
 we consider the curves that are extremal with respect to this functional.
First, we prove the variational formula, which is a nonlinear
differential equation relative to the general-affine curvature,  then
discuss some special solutions with reference to Chazy equations.  Second,
we compute the variational formula for a more generally defined
curvature functional.

\subsection{Extremal plane curves relative to the length functional}

 We have shown that there exists a unique frame $e = \{x, E_1, E_2\}$
 such that \eqref{pfaff} holds. Recall that  $\Omega$ denotes
the $3 \times 2$ matrix in \eqref{pfaff}.
Let $x_{\eta} (t)$ be a family of curves parametrized by $\eta$ 
around $\eta =0$ and $x_{0} = x$.
 We assume that $x_{\eta} (t) = x(t)$ outside a compact set $C$ 
 and $x_0(t)$ is parametrized by general-affine arc length.
For simplicity, we further assume that 
$x_{\eta}$ does not have an affine inflection point, {\it i.e.},
the invariant $\omega_s$ does not 
vanish anywhere for all $\eta$.
The length functional $L$ is given by 
\[
 L (\eta) = \int_C \omega_s (\eta).
\]
 For the sake of brevity, we use the notation $\delta$ to denote 
 the derivative with respect to $\eta$ evaluated at $\eta =0$:
\[
 \delta a = \frac{d a (\eta)}{d \eta}\Big|_{\eta =0}.
\]
 Then the curve $x$ is called  \textsl{general-affine extremal} if 
\[
 \delta L =0
 \]
 for any compactly supported deformation of $x$.

 We want to derive a differential equation for affine extremality.
 Since $\{E_1, E_2\}$ are linearly independent, there exists a 
 $3 \times 2$-matrix $\tau$ such that 
\[
 \delta \begin{pmatrix}x \\ E_1 \\E_2 \end{pmatrix} = 
 \tau \begin{pmatrix}E_1 \\E_2 \end{pmatrix}, \qquad
\tau = \begin{pmatrix} \tau_0^1 & \tau_0^2 \\
\tau_1^1 & \tau_1^2 \\ \tau_2^1 & \tau_2^2 
\end{pmatrix}
\]
 holds. 
Components of $\Omega$ and $\tau$ are denoted by $\omega_{\alpha}^{\beta}$
and $\tau_{\alpha}^{\beta}$, where $\alpha=0, 1, 2$ and  $\beta =1,2$.
Since $\delta d e = d \delta e$ with $e = \{x, E_1, E_2\}$, we have 
\[
 \delta \omega_{\alpha}^{\beta} -d \tau_{\alpha}^\beta
 =\sum_{\gamma=1,2}\tau_{\alpha}^{\gamma} \omega_{\gamma}^{\beta}
- \omega_{\alpha}^{\gamma} \tau_{\gamma}^{\beta}.
\]
 In terms of 
 entries of $\Omega$ and $\tau$, we have
 \begin{align}
\delta \omega_s - d \tau_{0}^1  &=  - \left(
 \frac{1}{2} k \tau_0^1 + \epsilon \tau_0^2 + \tau_1^1 \right)  \omega_s, \label{eq:0-1}\\
- d \tau_{0}^2  &=  \left(
 \tau_0^1 - \tau_1^2  -k \tau_0^2 \right)  \omega_s, \label{eq:0-2}\\
- \frac{1}{2}\delta (k\omega_s) - d \tau_{1}^1  &= - \left( \epsilon \tau_1^2 +\tau_2^1 \right) 
 \omega_s, \label{eq:1-1}\\
\delta \omega_s - d \tau_{1}^2  &=  \left(
 \tau_1^1 -\frac{1}{2} k \tau_1^2 -\tau_2^2\right)  \omega_s, \label{eq:1-2}\\
- \epsilon \delta \omega_s - d \tau_{2}^1  &=  
 \epsilon \left(\tau_1^1 + \frac{\epsilon}{2} k \tau_2^1 -\tau_2^2\right)  \omega_s, \label{eq:2-1}\\
- \delta (k \omega_s) - d \tau_{2}^2  &=  \left( \tau_2^1 + \epsilon \tau_1^2 \right)  \omega_s.\label{eq:2-2}
\end{align}
 Here we use $\omega_0^1 = \omega_1^2 =-\epsilon \omega_2^1= \omega_s$, 
 $\omega_0^2 =0$,  $\omega_1^1 = - 1/2 k \omega_s$,  and $\omega_2^2 = 
 - k \omega_s$.  
Then 
 adding \eqref{eq:1-2} and $- \epsilon$\eqref{eq:2-1}, 
\[
 2 \delta \omega_s -d \tau_1^2 + \epsilon d \tau_2^1 =  
 - \frac{1}{2} k (\tau_1^2 + \epsilon \tau_2^1) \omega_s
\]
 holds. Then adding  \eqref{eq:2-2} and $-2$\eqref{eq:1-1}, 
 we have 
\begin{equation}\label{eq:tau1221}
 2 d \tau_1^1 - d \tau_2^2 = 
3\epsilon ( \tau_1^2  + \epsilon \tau_2^1 ) \omega_s.
\end{equation}
Combining these equations, we have 
\begin{equation}\label{eq:deltaomegas}
  2 \delta \omega_s = d \tau_1^2 - \epsilon d \tau_2^1 + \frac{\epsilon}{6} k (-2 d \tau_1^1 + d \tau_2^2).
\end{equation}
 Recall that the deformation is compactly supported. 
 Then by using Stokes' theorem and integration by parts, 
 we have 
\begin{equation*}
  \delta L = -\frac{\epsilon}{12} \int_{C}  
 \left(-2\tau_1^1 + \tau_2^2\right) d k.
\end{equation*}
 We now compute $-2 \tau_1^1+\tau_2^2$ as follows. 
 From \eqref{eq:0-1} and \eqref{eq:1-2}, we have 
\begin{equation*}
d \tau_0^1- d \tau_1^2  = \left(2 \tau_1^1 - \tau_2^2  +
\epsilon \tau_0^2\right) \omega_s + \frac{1}{2}k \left(\tau_0^1  - \tau_1^2 \right) 
\omega_s.
\end{equation*}
 Inserting \eqref{eq:0-2} into the above equation, we have 
\begin{equation} \label{tau11tau22}
-2 \tau_1^1 + \tau_2^2 = {\tau_0^2}''  - \frac{3}{2}k {\tau_0^2}' 
+ \left(\epsilon +\frac{1}{2}k^2 - k'\right)\tau_0^2.
\end{equation}
 Here $\{'\}$ denotes the derivative with respect to the general-affine arc length 
 $\omega_s$, {\it i.e.} $a'= d a/\omega_s$ for a function $a$.
 Finally, using integration by parts again, 
\begin{equation}
  \delta L = -\frac{\epsilon}{12} \int_{C}  
 \left(k''' + \frac{3}{2}k k''  + \frac{1}{2} {k'}^2 + \frac{1}{2}  k^2 k' 
 + \epsilon k'\right) \tau_0^2 \omega_s
\end{equation} 
 holds. 
If we now take 
$x(t)+\eta \left\{v^1(t,\eta)E_1(t)+v^2(t,\eta)E_2(t)\right\}$
for the family of curves $x_\eta(t)$, where $v^1$ and $v^2$ are
arbitrary smooth functions with compact support
relative to $t$, then $\delta x =v^1(t,0)E_1 + v^2(t,0)E_2$. This
means that $\tau_0^2=v^2(t,0)$ is arbitrary for this family and the
vanishing of the integral implies the following proposition.

\begin{proposition}[\cite{Mihai2}]  \label{planeextremalprop}
 A {nondegenerate} plane curve without affine inflection point 
is general-affine extremal 
 relative to the length functional if and only if 
\begin{equation}\label{eq:variationformula}
 k''' + \frac{3}{2}  k k''  + \frac{1}{2} {k'}^2 + \frac{1}{2}  k^2 k' + 
 \epsilon k' =0
\end{equation}
 holds. 
In particular, any curves of constant general-affine curvature
are extremal.
\end{proposition}
We remark here that the differential equation
$(\ref{eq:variationformula})$ 
was first given in \cite[(33)]{Mihai2}, though some modifications
are necessary. The formula for $\epsilon=1$
was then rediscovered by S. Verpoort in \cite[p.432]{Ve} by 
making use of his general variational formula of equiaffine 
invariants: The differential equation is written in terms of 
both of the equiaffine curvature and the general-affine curvature,
and looks much simpler than $(\ref{eq:variationformula})$. As a result,
he proved the following corollary, which we now state in our setting.

\begin{corollary}[\cite{Ve}]\label{planeextremal}
Let $x(t)=(x_1(t), x_2(t))$ be a curve parametrized
by general-affine parameter $t$. Assume that it is general-affine 
extremal. Then, there exist constants $c_1$, $c_2$ 
and $c_3$ such that the general-affine curvature $k$ can be 
 written as
\[k= c_1x_1+c_2x_2 +c_3.\]
\end{corollary}

\noindent Proof. Let us consider the ordinary linear 
 differential equation
\[z'''+{3\over 2}kz'' + \left(\epsilon + {1\over 2}k'+{1\over 2}k^2\right)z'=0,
\]
with unknown function $z$. We note that this is nothing but of the
same form as the equation $(\ref{kode})$, therefore, $x_1$, $x_2$ are 
solutions. Also any constant is obviously a solution. On the other hand,
if $x$ is extremal, then the equation $(\ref{eq:variationformula})$ 
shows that $k(t)$ itself is a solution. Therefore, $k$ can be 
 expressed as claimed.

\medskip

We also have the following property on the curvature integral.

\begin{proposition}
 The variation of total curvature 
 on any compact interval always vanishes, {\it i.e.,}
\[
 \delta  \int_{C} k\omega_s =0
\]
 holds.
\end{proposition}

\noindent Proof. Adding \eqref{eq:1-1} and \eqref{eq:2-2}, 
\[
 - \frac{3}{2} \delta (k \omega_s)- d \tau_1^1 - d \tau_2^2 = 0
\]
 holds. Then Stokes' theorem implies the proposition.

\begin{remark}
{\rm 
 The differential equation \eqref{eq:variationformula} associated to 
 general-affine extremal plane curves is the third order nonlinear 
 differential equation of Chazy type. 
 It is known that 
 J. Chazy \cite{Ch} classified third order nonlinear ordinary 
 differential 
 equations of 
 Painlev\'e type, {\it i.e.} the solutions only admit poles as movable singularities.
 Then Chazy equations are classified into \textrm{I} to \textrm{X\!I\!I\!I} classes of 
 equations and the full list of equations can be found in \cite{Ba}.
 We here cite Chazy equations for 
 \textrm{I\!V, V} and \textrm{V\!I}, which are respectively given by
\begin{align}
& k''' +3k k''+ 3 {k'}^2  + 3k^2 k' - Sk' - S' k - T=0, \\
& k''' +2k k''+ 4 {k'}^2 +  2k^2 k' - 2R k' - R' k =0, \\
& k''' +k k''+ 5 {k'}^2 +  k^2 k' - 3Q k' - Q' k + Q''=0, 
\end{align}
 where $S, T, R, Q$ are certain analytic functions of $t$, \cite{Ba}.
 Then  \eqref{eq:variationformula} is clearly 
 of the form of the above Chazy equations with 
 the coefficients of $k k''$,  ${k'}^2$ and  $k^2 k'$ 
 replaced by the half integers $3/2$, $1/2$ and $1/2$, respectively, and 
with $S$, $T$, $R$ or $Q$ chosen properly. }
\end{remark}

\begin{example} \label{planeextremalexample}
{\rm
 We can see that the following $k(t)$ are solutions 
of \eqref{eq:variationformula}:
 \[
  k (t) = 3 \sqrt{2}\tanh (\sqrt{2}(t-c)) \quad\mbox{and}\quad
  k(t) = 3 \sqrt{2}\coth (\sqrt{2}(t-c))
 \]
for $\epsilon = 1$ (\cite[Example 14]{Ve}), and 
\[
k(t)=-3\sqrt{2}\tan(\sqrt{2}(t-c)),\quad 
k(t)=3\sqrt{2}\cot(\sqrt{2}(t-c))\quad \mbox{and}\quad
k(t)=\pm\sqrt{2}+\frac{3}{t-c}
\]
for $\epsilon =-1$. For each of these solutions, we can
compute the associated plane curve by integrating the differential
equation, using computer software.  Since the expression is not
simple, we give here one example for the case $\epsilon=-1$ and
$k(t)=\sqrt{2}+3/t$, and the curve is written as $(x_1,x_2)$ for $t>0$:
\begin{align*}
x_1 &= 3\sqrt{2} -{1\over t}, \\
x_2 &= {\sqrt{3\pi}(\sqrt{2}-6t)\ 
{\rm erfi}\left({\sqrt{3t}\over 2^{1/4}}\right)\over t}
+ {62^{1/4}\exp\left({3 t \over \sqrt{2}}\right)\over \sqrt{t}},
\end{align*}
where erfi is the error function defined by
\[
 {\rm erfi}(x) = {2 \over \sqrt{\pi} } \int_0^{x} e^{t^2} \ dt.
\]
}
\end{example}

\subsection{Extremal problem for a generalized curvature functional}
 More generally, one can consider a variational problem for the following 
 curvature functional:
\begin{equation}\label{eq:generalfunc}
 F(\eta) = \int_C f(k) \omega_s,
\end{equation}
 where $f$ is a smooth function of one variable and
$\eta$ is the parameter for the variation of curves. Then we have
\begin{equation}\label{eq:generalfunc2}
 \delta F = \int_C f^{\prime}(k) (\delta k) \omega_s + 
 f(k) \delta \omega_s.
\end{equation}
 The computation of $\delta k $ is done as follows:
 Adding \eqref{eq:1-2} and $\epsilon$\eqref{eq:2-1}, we have 
 \[
  - d \tau_1^2 - \epsilon d \tau_2^1 = 
 2 (\tau_1^1 - \tau_2^2) \omega_s - \frac{1}{2} k (\tau_1^2 - \epsilon \tau_2^1)
 \omega_s.
 \]
 By taking a derivative of this equation and 
by using \eqref{eq:tau1221}, we have
\begin{equation}\label{eq:ddot}
 \frac{\epsilon}{3}\left(-2 {\tau_1^1} + {\tau_2^2}\right)'''=
  2 ({\tau_1^1} - {\tau_2^2})' - \frac{1}{2} k' ( \tau_1^2 - \epsilon \tau_2^1)
 - \frac{1}{2}k ({\tau_1^2} - \epsilon {\tau_2^1})'.
\end{equation}
Subtracting \eqref{eq:2-2} from \eqref{eq:1-1}, we have 
\[
 \frac{1}{2} (\delta k) \omega_s + \frac{1}{2} k \delta \omega_s - d \tau_1^1 + d \tau_2^2
= -2 (\tau_2^1 + \epsilon \tau_1^2) \omega_s.
\]
Then, using \eqref{eq:deltaomegas}
we see that
\begin{equation*}
 \delta k =  -k \left( \frac{1}{2} \left(\tau_1^2 -\epsilon \tau_2^1\right)' 
 + \frac{\epsilon}{12} k (-2 \tau_1^1 + \tau_2^2)'\right)
+ 2 (\tau_1^1 - \tau_2^2)' -4 (\tau_2^1 + \epsilon \tau_1^2).
\end{equation*}
Again, by using \eqref{eq:tau1221},  
\begin{align*}
 \delta k = - \frac{1}{2} k \left( \tau_1^2 -\epsilon \tau_2^1\right)'  
 + \left(\frac{4}{3} - \frac{\epsilon}{12}k^2\right)
\left( -2 {\tau_1^1} + {\tau_2^2}\right)' + 2 ( \tau_1^1 - \tau_2^2)',
\end{align*}
and finally, by using \eqref{eq:ddot}, we have 
\begin{equation}\label{eq:deltak}
 \delta k =  \frac{1}{2}k'
 \left(\tau_1^2 - \epsilon \tau_2^1\right)+
 \frac{\epsilon}{3}\left(-2 {\tau_1^1} + {\tau_2^2} \right)'''
+
 \left(\frac{4}{3}- \frac{\epsilon }{12}k^2\right) \left( -2 {\tau_1^1} + {\tau_2^2}\right)'. 
\end{equation}
Then, by inserting the expressions $\delta\om_s$ 
\eqref{eq:deltaomegas} 
and $\delta k$ \eqref{eq:deltak} into \eqref{eq:generalfunc2}, and
by using integration by parts, we get
\begin{align*}
\delta F  = 
\int_C 
 \dot f (k)\left\{\frac{\epsilon}{3}\left(-2{\tau_1^1} + {\tau_2^2} \right)'''\right.
  +  \left(- \frac{\epsilon k^2}{12} + \frac{4}{3}\right) &\left. \left( -2 {\tau_1^1} 
 + {\tau_2^2}\right)'\right\} \omega_s
 \\ 
& 
 + f(k) \left\{\frac{\epsilon}{12} k \left(-2 {\tau_1^1} + {\tau_2^2}\right)'\right\} 
 \omega_s.
\end{align*}
 Again applying integration by parts, we have 
\[
 \delta F = 
-\frac{\epsilon}{12}\int_C G (-2 \tau_1^1 + \tau_2^2) \omega_s,
\]
where 
\begin{equation}\label{eq:G}
 G = 4 \ddddot f(k){k'}^3 + 12 \dddot f(k) k' k'' + 
\ddot f(k) (4 k''' - k' k^2 + 16\epsilon k'  ) 
 - \dot f (k) k k' + f(k)  k'.
\end{equation}
Thus, by use of $(\ref{tau11tau22})$, we have the following theorem.
\begin{theorem}  \label{genvarform}
  A plane curve without affine inflection points is 
general-affine extremal with respect to 
 the curvature functional \eqref{eq:generalfunc} if and only if 
 \begin{equation} \label{Gextremal}
  G'' + \frac{3}{2}{G'} k + \frac{1}{2}G k'  + \frac{1}{2}G k^2+
 \epsilon G =0
 \end{equation}
 holds, where $G$ is the function defined in \eqref{eq:G}.
\end{theorem}

\begin{remark} 
Variation of energy integral. 
{\rm
When $f={1\over 2}k^2$, the integral $F$ may be called the energy integral.
For this $f$, we see that 
\[ G=4k'''-{3\over 2}k^2k' + 16\epsilon k'\]
and the equation $(\ref{Gextremal})$ give an extremal curve relative to
the energy functional.
}
\end{remark}

\section{How to find plane curves with given general-affine curvature}
\label{sec:graphimmersion}
In Section \ref{subsc:graph}, we have derived 
the expression \eqref{curv} of the general-affine curvature
for a graph immersion $x(t) = (t, f(t))$
 with $\mu = (f^{\prime \prime})^{-2/3}>0$.
Making use of this expression, we study how to find 
a graph immersion of plane curves 
with given general-affine curvature,
 by considering the following nonlinear differential equation directly, 
\begin{equation} \label{muconst}
 \mu (\mu''')^2 =  -\epsilon {k^2\over 2}(\mu'')^3.
\end{equation}
 We regard the function $\mu'$ of $t$ as a function
of $\mu$ and set
\[ w(\mu) = \mu'(t) ={d\mu\over dt}.\]
Then, by the chain rule, we have
\[ \mu'' = w\dot{w},\quad \mu'''=w\dot{w}^2+ w^2\ddot{w}.\]
 Hence, the equation $(\ref{muconst})$ is written as
\begin{equation} \label{wode}
  \mu w^2(\dot{w}^2+w\ddot{w})^2 + \epsilon {k^2\over 2}w^3\dot{w}^3=0,
\end{equation}
which can be reduced to the Abel equation as follows:

\noindent {(i) First reduction}:
We introduce $s$ by setting
\[ w(x) = \pm \exp \left( -\epsilon \int s^2 dx\right).\]
Here we choose the sign properly, depending on the function $w$.
Then, we get the equation
\begin{equation} \label{mug}
  8x(-\epsilon \dot{s} + s^3)^2 - k^2s^4 =0, \quad \quad (x>0).
\end{equation}
Therefore, the original differential equation $(\ref{muconst})$ is
equivalent to
\begin{equation}\label{eq:abel1}
  \epsilon \dot s={k\over 2\sqrt{2x }}s^{2} + s^3,
\end{equation}
which is an \textsl{Abel equation of the first kind}.
%
%

 It is easy to see that for constant $k<0$ with $\epsilon =- 1$
 or $k\leq -4$ with $\epsilon = 1$, the solution $s$ 
 of \eqref{eq:abel1} can be explicitly obtained as
\[
 s(x) = \frac{a}{\sqrt{2 x}} \quad \mbox{with}\quad  a 
= \frac{-k \pm\sqrt{- 16\epsilon +k^2}}{4}.
\]
 The corresponding {curves} are given in Examples \ref{exp} and \ref{power}.
 Moreover, in the case of $k =0$ (both $\epsilon = \pm 1$),  
 the solution $s$ can be obtained as
\[
 s(x) = \frac{1}{\sqrt{\epsilon(a -2 x)}},
\]
 where $a$ is some constant.
 The corresponding {curves} are given in Example \ref{ellzero}.
 On the contrary, in the case of $-4 < k < 0$ for $\epsilon =1$, 
 the solution $s$ of \eqref{eq:abel1} is not easy to write down explicitly. 
 The corresponding {curves} are given in Example \ref{logsp}.

\medskip

\noindent {(ii) Second reduction}: We define $s$ by 
\[w(x)= \pm \exp\left(-\epsilon \int s^{-2}dx\right),\]
by choosing the sign properly.
Then a straightforward computation shows that the equation $(\ref{wode})$
is transformed to
\[ -k^2 s^2 + 8x(\epsilon +s\dot s)^2=0,\]
which is equivalent to 
\begin{equation}\label{eq:abel2}
s\dot s = \frac{k}{2 \sqrt{2 x}} s - \epsilon.
\end{equation}
This is a particular case of the \textsl{Abel equation of the second kind}.
We refer to \cite[Section 1.3.2]{PZ} for integrable Abel equations.

\begin{theorem} \label{abeleq}
 For any general-affine plane curve with graph immersion $(t, f(t))$,
 there exists a function $s$ given as above such that
 $s$ satisfies the Abel equation of the first kind or second kind,
 \eqref{eq:abel1} or \eqref{eq:abel2}, 
 respectively. Conversely, for given any function $k$, a
 solution $s$ of \eqref{eq:abel1} or \eqref{eq:abel2} 
 gives rise to a plane curve of graph immersion 
 $(t, f(t))$ with general-affine curvature $k$.
\end{theorem}


\section{General-affine curvature of space curves}

The equiaffine treatment of space curves as well as the projective
treatment of space curves are classically known. However, it seems that
a general-affine treatment of space curves is not fully developed.
In this section, we will introduce several notions
such as curvature, length parameter and ordinary differential equation
associated with space curves from a general-affine point of view.

\subsection{Choice of frames for space curves and general-affine curvatures}
\label{subsc:3frames}

Let $x:t \longrightarrow x(t)\in {\bf A}^{3}$ be a curve
in a $3$-dimensional affine space with parameter $t$
and let $e=\{e_1,e_2,e_3\}$ be a frame along $x$; 
it is a set of independent vectors of ${\bf A}^{3}$.
The vector-valued $1$-form $dx$ is written as
\begin{equation} \label{dx3}
dx = \om^1 e_1 + \om^2 e_2 + \om^3e_3,
\end{equation}
and the dependence of $e_i$ on the parameter is described by the equation
\begin{equation} \label{dei3}
de_i = \sum_j \om_i^j e_j,
\end{equation}
where $\om^j$ and $\om_i^j$ are $1$-forms as before in the 2-dimensional
case and $1\le i,j \le 3$. We call $\{\om^i, \om_i^j\}$ the coframe.

We assume in the following that the curve is nondegenerate in the 
sense that the vectors $x'$, $x''$ and $x'''$ are linearly independent
and that $\om^2=\om^3=0$ and $\om_1^3=0$, so that $e_1$ is 
tangent to the curve and that $\{e_1, e_2\}$ is the first osculating
space of the curve. We write $\om^1 = \om$ for simplicity.

Let $\tilde{e}=\{\tilde{e}_1, \tilde{e}_2, \tilde{e}_3\}$ be 
another choice of such a frame. Then, it can be written as
\[\tilde{e}_1 = \lam e_1,\qquad \tilde{e}_2 = \mu e_1 + \nu e_2,\qquad
\tilde{e}_3=\alpha e_1+ \beta e_2+\gamma e_3,\]
where $\lam\nu\gamma \neq 0$.
The associated coframe is written as $\tilde{\om}$ and $\tilde{\om}_i^j$, 
which satisfies
\[ dx=\tilde{\om}\tilde{e}_1, \qquad d\tilde{e}_i = \sum_j \tilde{\om}_i^j\tilde{e}_j.
\]
Then we have
\begin{equation} \label{shiki0}
\tilde{\om} = \lam^{-1}\om.
\end{equation}
Since $d\tilde{e}_1$ is represented in two ways, one being
\[ d\tilde{e}_1 = (d\lam) e_1+ \lam (\om_1^1e_1+\om_1^2e_2)\]
and the other being
\[ d\tilde{e}_1 = \tilde{\om}_1^1(\lam e_1)+\tilde{\om}_1^2(\mu e_1+\nu e_2),\]
by comparing the coefficients of $e_1$ and $e_2$ in these expressions, we get
\begin{eqnarray} 
& \nu \tilde{\om}_1^2= \lam \om_1^2, & \label{shiki1} \\
& \lam \tilde{\om}_1^1 + \mu \tilde{\om}_1^2 = d\lam  + \lam \om_1^1. &
\label{shiki2}
\end{eqnarray}
Similarly, by considering $d\tilde{e}_2$, we have
\def\tom{\tilde{\omega}}
\begin{eqnarray}
&\gamma\tom_2^3 = \nu \om_2^3,& \label{shiki3} \\
& \nu \tom_2^2+\beta \tom_2^3 = d\nu + \mu\om_1^2 + \nu \om_2^2,& \label{shiki4}
\\
& \lam \tom_2^1 + \mu\tom_2^2 + \alpha \tom_2^3 =
  d\mu + \mu\om_1^1 + \nu\om_2^1,& \label{shiki5}
\end{eqnarray}
and by $d\tilde{e}_3$ we have 
\begin{eqnarray} 
& \gamma\tom_3^3 = d\gamma + \beta\om_2^3+\gamma\om_3^3, & \label{shiki6} \\
&\nu\tom_3^2+\beta\tom_3^3=d\beta+\alpha\om_1^2+\beta\om_2^2+\gamma\om_3^2, &
\label{shiki7} \\
&\lambda\tom_3^1+\mu\tom_3^2+\alpha\tom_3^3 =d\alpha + \alpha\om_1^1
  + \beta\om_2^1 + \gamma\om_3^1.& \label{shiki8}
\end{eqnarray}
First note that, from the generality assumption, we have $\om_1^2\neq 0$
and $\om_2^3\neq 0$. Then, by an appropriate choice of $\nu$ and $\gamma$,
in view of $(\ref{shiki1})$ and $(\ref{shiki3})$, we can assume 
that $\tom_1^2=\tom$ and $\tom_2^3=\tom$. Hence, we
can restrict our consideration to the case 
\[\om_1^2=\om\quad {\rm and} \quad \om_2^3=\om\]
in the following. In particular,
\begin{equation} \label{unimo}
\nu=\lambda^2\quad {\rm and}\quad \gamma=\lambda^3
\end{equation}
are necessary. We next see that, from $(\ref{shiki2})$,
$(\ref{shiki4})$ and $(\ref{shiki6})$, we have
\begin{align*}
 2 \tom_1^1 - \tom_2^2 &= 2 \om_1^1 - \om_2^2 - 3 \lam^{-2} \mu \om + \lam^{-3} \beta \om, \\
 3\tom_1^1 - \tom_3^3 &= 3 \om_1^1 - \om_3^3 - 3 \lam^{-2} \mu \om - \lam^{-3} \beta \om. 
\end{align*}
Thus an appropriate choice
of the parameters $\mu$ and $\beta$ makes the identities
$\tom_2^2=3 \tom_1^1$ and $\tom_3^3=2 \tom_1^1$ hold. 
To keep this condition it is necessary to have $\mu=\beta=0$. Now
 \eqref{shiki5} can be rephrased as
\[
  \lam \tom_2^1 + \alpha \tom_2^3 = \lam^2 \om_2^1,
\]
 and we choose $\alpha$ so that $\tom_2^1 =0$. Thus, we can assume that 
 $\om_2^1=0$ and $\alpha=0$ in the following.
 Moreover \eqref{shiki2} is 
\[
 \tom_1^1 = \lam^{-1} d \lam  + \om_1^1,
\]
and we can choose $\lam$ so that $\tom_1^1=0$. 
Therefore $\om_1^1=0$, and to keep this condition, 
 $\lam$ is a non-zero constant.
With these considerations, the last identities $(\ref{shiki7})$
and $(\ref{shiki8})$ turn out to be
\[\tom_3^2 = \lambda\om_3^2 \quad {\rm and}\quad
\tom_3^1=\lambda^2\om_3^1,\]
respectively. We set
\begin{equation} \label{spcinv}
\om_3^2 = -\ell \om, \quad \om_3^1=-m\om,
\end{equation}
and similarly for $\tom_3^2$ and $\tom_3^1$. Then, we have the
covariance
\begin{equation} \label{spccov}
\tilde{\ell}=\lambda^2\ell, \quad {\rm and}\quad
\tilde{m} = \lambda^3 m.
\end{equation}
Thus, we have seen that, given a nondegenerate curve $x$, there 
 exists a frame $e$ with the coframe of the form
\begin{equation}
\label{eq:coframeeqaff}
 \begin{pmatrix}
\omega & 0 &0  \\ 
0 & \omega & 0 \\ 
0& 0 &  \omega\\
-m \omega & - \ell \omega & 0 
 \end{pmatrix}.
\end{equation}

We remark here that, in the equiaffine treatment of space curves, 
the scalars $\ell$ and $m$ above are known to be absolute invariants,
called the \textsl{equiaffine curvature} and the \textsl{equiaffine torsion}, 
respectively; we refer 
to Section \ref{3equiaffine}.
In this paper, we call the point where $\ell=0$ an \textsl{affine
inflection point}.

In the following we assume $\ell\neq 0$ 
and 
let $\epsilon$ denote the sign of $\ell$: 
\[\epsilon={\rm sign}(\ell).\]
It is an invariant of the curve.
Then we define the \textsl{general-affine length element} by
\begin{equation}  \label{galength}
\om_s = \sqrt{\epsilon \ell}\om,
\end{equation}
which is well-defined independent of the frame in view of $(\ref{spccov})$,
and a parameter $s$ for which
$ds = \om_s$ holds is the \textsl{general-affine length parameter}
determined up to an additive constant.

\begin{definition} \label{spacecurv}
We call the scalar function $k$ defined by
\[ {d\ell\over \ell} = k\om_s\]
the \textsl{first general-affine curvature}.
In other words,
\begin{equation}\label{kforspacecurve} 
k={d\log\ell\over ds}.
\end{equation}
We call the scalar function $M$ defined by 
\begin{equation} \label{gasecond}  
M={m\over (\epsilon\ell)^{3/2}}
\end{equation}
the \textsl{second general-affine curvature} of the space curve.
\end{definition}
Both curvatures defined above are absolute invariants.

We next define a new frame $\{E_1, E_2, E_3\}$ by setting
\[E_1 = {1\over (\epsilon\ell)^{1/2}}e_1,\qquad E_2={1\over \epsilon \ell}e_2,
\qquad E_3={1\over (\epsilon\ell)^{3/2}}e_3.\]
It is easy to see that this frame does not depend on the choice of $\lambda$;
hence, it is determined uniquely.

Thus we have proved the following:

\begin{proposition} Assume $\ell\neq 0$. Then, the frame $\{E_1, E_2, E_3\}$ 
is uniquely defined from the immersion and it satisfies the Pfaffian equation
\begin{equation} \label{pfaff3}
 d \left(\begin{array}{c} x \\ E_1 \\ E_2 \\ E_3\end{array}\right)
=\Omega
\left(\begin{array}{c} E_1 \\ E_2 \\E_3 \end{array}\right),
\qquad
\Omega =\left(
\begin{array}{ccc}
 \om_s    &   0  &  0 \\
-{1\over 2}k\om_s & \om_s & 0 \\
0 & -k\om_s  & \om_s \\
- M \om_s & -\epsilon \om_s & -{3\over 2}k\om_s
\end{array}\right),
\end{equation}
where $\om_s$ is the general-affine length form,
$k$ and $M$ are the first and second general-affine curvatures,
respectively, and $\epsilon={\rm sign}(\ell)$. 
\end{proposition}

By use of this choice of frame, we can see the following lemma,
by a similar reasoning to that for Lemma \ref{lem:kode}.
\begin{lemma} \label{gaspclemma}
The immersion $x$ satisfies the ordinary differential
equation
\begin{equation} \label{gaspcode}  
x''''+3kx''' +\left(2k'+{11\over 4}k^2+\epsilon\right)x'' 
  + \left(M+{1\over 2}\epsilon k + {1\over 2}k''
          +{7\over 4}kk' + {3\over 4}k^3\right)x'=0,
\end{equation}
relative to a general-affine length parameter.
\end{lemma}

In the definition of the curvature, 
we had an ambiguity of orientation of the chosen parameter:
by the change of the parameter from $t$ to $-t$, the equation transforms 
 to
\[ 
x''''-3kx''' +\left(-2k'+{11\over 4}k^2+\epsilon\right)x'' 
  + \left(-M-{1\over 2}\epsilon k - {1\over 2}k''
          +{7\over 4}kk' - {3\over 4}k^3\right)x'=0.
\]
Namely, the transform $(k,M)\rightarrow (-k,-M)$ keeps the form
of the equation.

Thus, up to this ambiguity, we have the following theorem.

\begin{theorem} \label{spacenatural}
Given functions $ k(t)$ and $M(t)$ of a parameter $t$, and $\epsilon=\pm 1$, 
there exists a nondegenerate space curve $x(t)$ for which $t$ is
a general-affine length parameter, $k$ is the first 
general-affine curvature, $M$ is the second general-affine
curvature, and $\epsilon$ is the sign of $\ell$, 
uniquely up to a general-affine transformation.
\end{theorem}

Analogously to the case of plane curves, we have the following property
on the total general-affine curvature:

\begin{corollary} \label{spacetotalcurv}
Assume that the curve $C$ is nondegenerate and closed, and 
has no affine inflection point.
Then, the total curvature
$\int_C k\om_s$ vanishes. In particular, such a curve has at least
two general-affine flat points.
\end{corollary}

\comment{
\begin{remark}
{\rm
We have chosen a certain frame 
and determined the coframe to be $(\ref{pfaff3})$
under the condition $\ell\neq 0$;
we call this process a normalization of frame. As is seen from the above,
however, the normalization is not unique. In fact, we can show that
there is another 
normalization by which the coframe has the form
\[
\left(
\begin{array}{ccc}
 \om_s    &   0  &  0 \\
-{1\over 2}k\om_s & \om_s & 0 \\
-\epsilon p \om_s & -k\om_s  & \om_s \\
- M \om_s & -\epsilon \om_s & -{3\over 2}k\om_s
\end{array}\right),
\]
 where $p$ is not equal to $-2$, which we do
 not represent here.
 Moreover, it seems that 
 Mih$\breve{\rm a}$ilescu 
 gave another normalization in \cite{Mihai3}
under the condition that the curve does not belong to a linear complex
instead of assuming $\ell\neq 0$. 
}
\end{remark}
}

\subsection{Computation of general-affine curvatures of space curves}

Let $t\longrightarrow x=x(t)\in {\bf A}^3$ be a nondegenerate curve 
such that the vectors $x'$, $x''$ and $x'''$ are linearly independent. 
Since $x''''$ is written as a linear combination of $x'$, $x''$ and $x'''$,
there are scalar functions $a=a(t)$, $b=b(t)$ and $c=c(t)$ such that
\begin{equation} \label{x3ode}
x'''' = a x''' + b x'' + cx'.
\end{equation}
We give a formula to compute these coefficients 
by use of the general-affine curvatures of such a curve.
The method is similar to that used for plane curves.

Since $dx=x'\, dt$, the frame vector $e_1$ is a scalar multiple of $x'$:
\begin{equation}\label{eq:3frame}
 dx = \om\, e_1;\qquad e_1=\lambda x',\quad \om = \lambda^{-1}dt.
\end{equation} 
Then, the differential
\[ de_1=(\lambda^2 x'' + \lam\lambda' x')\,  \om\]
implies  that the second frame vector is
\[ e_2=\lambda^2x'' + \lam\lambda' x'.\]
The derivation of $e_2$ is
\[
de_2=(\lambda^3x'''+3\lam^2\lam'x''+(\lam^2\lam''+\lam{\lam'}^2)x')\om,\]
which is equal to $\om e_3$:
\[e_3 =(\lambda^3x'''+3\lam^2\lam'x''+(\lam^2\lam''+\lam{\lam'}^2)x').\]
Its derivation is
\[de_3=\left((\lambda^3a+6\lam^2\lam')x'''+(\lam^3b+7\lam{\lam'}^2
+4\lam^2\lam'')x''
+(\lam^3c+4\lam\lam'\lam''+\lam^2\lam'''+{\lam'}^3)x'
\right)dt\]
by use of $(\ref{x3ode})$. Since $de_3$ has no $e_3$-component, we have
\begin{equation} \label{eq:3lambda}
\lambda a + 6\lambda'=0, \quad {\it i.e.}\quad
   \lambda = e^{-{1\over 6}\int a(t) dt}
\end{equation}
up to a multiplicative constant. Then, $de_3$ is written as
\[ de_3 = (\lam^2b+7{\lam'}^2+4\lam\lam'')\om e_2 +
(\lam^3c - \lam^2\lam'b-6{\lam'}^3+\lam^2\lam''')\om e_1.
\]
By the definition in \eqref{spcinv}, we have
\begin{equation} \label{3ell} 
\ell = -(\lam^2b + 7{\lam'}^2 + 4\lam \lam'').
\end{equation}
Also, by the definition of $m$, we have
\begin{equation} \label{3m}
m = - \lam^3 c+\lam^2\lam'b+6{\lam'}^3 - \lam^2\lam'''.
\end{equation}
We now assume that $\ell\neq 0$ and recall that
$\epsilon={\rm sign}(\ell)$.
Then, we have
\[
ds^2 = \epsilon\ell\om \om
=  -\epsilon \left(b+7{\lambda'^2\over \lambda^2}+4\frac{\lam''}{\lam}\right) dt^2.
\]
In terms of $a$ and $b$, 
\begin{equation}\label{dstwo3}
ds^2 = -\epsilon \left(b+{11\over 36}a^2 - {2\over 3}a' \right) dt^2.
\end{equation}
Hence, a length parameter $s$ which is a function of $t$ is
obtained by solving the equation
 \[ \left({ds \over dt}\right)^2 
   = -\epsilon \left(b+{11\over 36}a^2 - {2\over 3}a'\right).\]
If, in particular, $t$ itself is a length parameter, then we must have
\begin{equation}\label{lengthspace}
\ell=\epsilon \lambda^2, \qquad 
 b=-\epsilon - {11\over 36}a^2 + {2\over 3}a'.
\end{equation}
By definition, the first curvature $k$ is
\begin{equation}\label{3curv}
k=-{1\over 3}a.
\end{equation}
We next treat the second curvature $M$ defined in 
{$(\ref{gasecond})$}:
From the formula $(\ref{3m})$ above,  
\begin{equation} \label{Mvalue}
M = -c +{\lam'\over \lam}b+6\left({\lam'\over \lam}\right)^3
- \frac{\lam'''}{\lam}.
\end{equation}
Hence, by $(\ref{eq:3lambda})$, we can see that 
\begin{equation} \label{cvalue}
c= -M + {1 \over 6}a\epsilon + {1\over 6}a''-{7\over 36}aa'+{1\over 36}a^3.
\end{equation}
Thus, we have seen that the differential 
equation $(\ref{x3ode})$ agrees with the equation $(\ref{gaspcode})$.


For another parameter $\sigma=\sigma(t)$, we write
\[ y(\sigma) = x(t).\]
Then, using the notation $\{\,\dot{}\,\}$ for the derivation by $\sigma$
and $\{\,{}'\,\}$ for the derivation by $t$, we see that
\begin{eqnarray*}
x' &=& \dot{y}\sigma', \\
x'' &=& \ddot{y}\sigma'^2 + \dot{y}\sigma'',\\
x''' &=& \threedot{y}\sigma'^3+3\ddot{y}\sigma'\sigma''+\dot{y}\sigma''',\\
x'''' &=& \fourdot{y}\sigma'^4 + 6\threedot{y}\sigma'^2\sigma''
+\ddot{y}(3{\sigma''}^2+4\sigma'\sigma''')+\dot{y}\sigma''''.
\end{eqnarray*}
Making use of these formulas, we can show that
\begin{equation} \label{yode}
\fourdot{y} = A(\sigma)\threedot{y}+B(\sigma)\ddot{y}+C(\sigma)\dot{y},
\end{equation}
where
\begin{eqnarray}
A(\sigma) &=& \left(a-6{\sigma''\over \sigma'}\right){1\over \sigma'}, 
\label{newspca}\\
B(\sigma) &=& \left( b+ 3a{\sigma''\over \sigma'}-3\left({\sigma''\over \sigma'}\right)^2-4{\sigma'''\over \sigma'}\right){1\over \sigma'^2}, 
\label{newspcb} \\
C(\sigma) &=& \left( c+b{\sigma''\over \sigma'} + a{\sigma'''\over \sigma'}
 - {\sigma''''\over \sigma'}\right){1\over \sigma'^3}.
\label{newspcc}
\end{eqnarray}

The differential polynomials that appeared in the representation
of $b$ and $c$ in \eqref{lengthspace} and \eqref{cvalue}
have a covariant property with respect to this change of parameters:

\begin{lemma} By the change of parameter, the following covariant 
relations hold.
\begin{equation} \label{ABCformula}
\begin{array}{rcl} \displaystyle
B-{2\over 3}\dot{A}+{11\over 36}A^2 &=&
\displaystyle
\left(b-{2\over 3}a' + {11\over 36}a^2\right){1\over \sigma'^2}, \\
\displaystyle
C-{1\over 6}\ddot{A}+{7\over 36}A\dot{A}-{1\over 36}A^3 &=&
\displaystyle
\left(c-{1\over 6}a''+{7\over 36}aa'-{1\over 36}a^3\right){1\over \sigma'^3}
\\
&&\qquad
\displaystyle
+ \left(b-{2\over 3}a' + {11\over 36}a^2\right){\sigma''\over \sigma'^4}.
\end{array}
\end{equation}
\end{lemma}

Thanks to the formulas above, we can compute curvatures
according to a procedure similar to that in Section \ref{subsectionga}.
\medskip

\begin{example} Viviani's curve. {\rm This curve is given by the mapping
\[(1+\cos(2t),\sin(2t),2\sin(t)).\]
The associated differential equation is
\[ x'''' = - \tan(t)x''' - 4x'' - 4\tan(t)x',\]
which is singular at $t$ with $\cos(t)=0$; in the left figure,
this corresponds to $z = \pm 2$.
A simple calculation shows the identity
\[-b-{11\over 36}a^2+{2\over 3}a'=
{5(31 \cos(t)^2-7)\over 36\cos(t)^2};\]
hence, at the values $t$ with $\cos(t)^2=7/31$, 
the general-affine length parameter cannot be defined, 
namely, $\ell=0$ at these values; in the left figure, there
correspond to the four points with $z=\pm 1.75...$. 
Except for these six values of $t$ (we marked these points as dots in 
 the figure), 
$\epsilon$ is determined and the curvatures are computable. 
The first curvature $k$ has the absolute value
\[ 
\frac{2 |\sin t| (49- 31 \cos^2 t)}{\sqrt{5}|31 \cos^2 t -7|^{3/2}}.
\]
}
\end{example}

\begin{example} Torus knot. {\rm The mapping 
\[x=((4+\cos(3t))\cos(t), (4+\cos(3t))\sin(t), \sin(3t))\]
defines one of the torus knots. The equation is computed as
\begin{eqnarray*}
 x'''' &=& {-3\sin(3t)(12T^2-152T-35)\over P}x'''
+{2(52T^3-178T^2+562T-891)\over P}x'' \\
&& +{12\sin(3t)(8T^2+82T+281)\over P}x',
\end{eqnarray*}
where
\[T=\cos(3t)\quad {\rm and}\quad P= 4T^3-76T^2 - 35T +198.\]
Since $P>0$ for all value $t$, the equation is non-singular.
With computer assistance, we can see that the length parameter is
well-defined and $\epsilon=1$, and that $k$
has period $2\pi/3$ and symmetry $k(t)=-k(2\pi/3-t)=-k(-t)$,
with values $-4<k<4$.
}\end{example}

\begin{figure}
\centering
\begin{tabular}{cc}
\includegraphics[width=3.8cm]{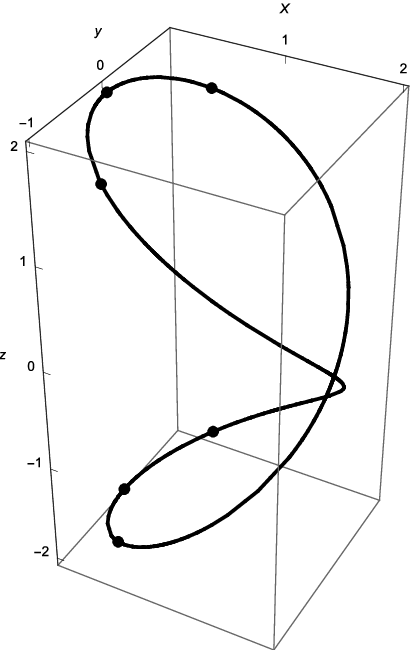} &
\includegraphics[width=5cm]{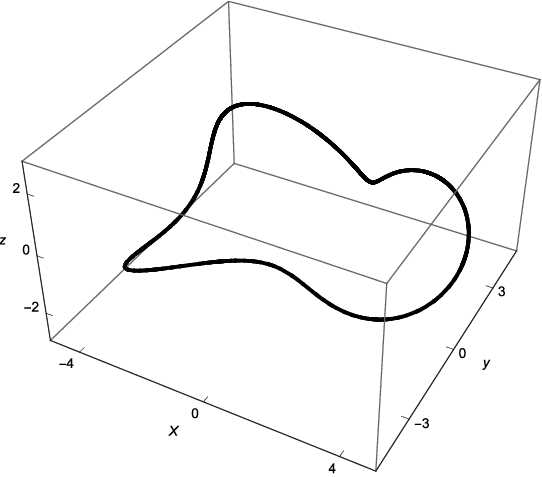} \\
{\rm Viviani's curve} & {\rm Torus knot}
\end{tabular}
\end{figure}


\subsection{Space curves with constant curvatures} \label{spcconst}

The curves with constant $k$ and $M$ have a special interest, because
such a curve is an orbit of a $1$-parameter subgroup of general-affine
motions. In \cite[p. 36--39]{Sc} a classification of such groups is given.
We here give a list of such curves by use of the differential equation
treated in the previous subsection.

If the curvatures are constant, then the differential equation
\[ x'''' = ax''' + bx'' + cx',\]
has constant coefficients $a$, $b$ and $c$.
Conversely, assume that a curve satisfies a differential equation with
constant coefficients. Then, 
the length parameter $s$ is obtained by the identity
\[ ds^2 = -\epsilon\left(b+{11\over 36}a^2\right) dt^2,\]
where
\[\epsilon ={\rm sign}\left(-b-{11\over 36}a^2\right).\]
We set
\[ q = \sqrt{-\epsilon \left(b+{11\over 36}a^2\right)}.\]
Then, the first curvature is
\[ k = -{1\over 3}\left( \frac{a}{q}\right),\]
and the second curvature is
\[ 
M=-{c\over q^3} + {\epsilon \over 6}\left({a\over q} \right)
  +{1\over 36}\left({a\over q}\right)^3.
\]
Therefore, any differential equation of the form above
with constant coefficients
defines a curve with constant general-affine curvatures. It is
sufficient to solve
\[ y'''=ay''+by'+cy\]
and integrate it to get $x$.
If the function $y=e^{\lam t}$ is a solution of the equation, then
$\lam$ is a root of the algebraic equation
\[\lam^3-a\lam^2-b\lam-c=0.\]
Depending on whether $\lam$ is a single or multiple root or a real
or imaginary root, the form of the solution varies. 
Without showing detailed computations, 
we list the result of the classification as follows,
which agrees with the classification of the $1$-parameter subgroup of 
general-affine motions.
First, we start with some examples.

\begin{example} \label{helix}
{\rm Curves with $k=M=0$. The equation is
\[x'''' = - \epsilon x''.\]
When $\epsilon=1$, the curve is $x=(t, \cos t, \sin t)$ called
a \textsl{circular helix} and, when $\epsilon=-1$, the curve is
$x=(t,\cosh t, \sinh t)$, called a \textsl{hyperbolic helix}.
}
\end{example}

\begin{example} \label{logspiral}
Logarithmic spiral. 
{\rm This curve is given by
 $x = (e^{-2\lam t}, e^{\lam t}\cos pt, e^{\lam t}\sin pt)$. 
 The equation is
\[ x'''' = (3\lam^2-p^2)x'' - 2\lam(p^2+\lam^2)x'.\]
 We see that $\epsilon=1$ (resp. $\epsilon=-1)$ when
$3\lam^2< p^2$ (resp. $3\lam^2> p^2$), $k=0$, 
and $M=2\lam(p^2+\lam^2)|3\lam^2-p^2|^{-3/2}$.
}
\end{example} 

\begin{example} {\rm Curves given by $x=(e^{\lambda t}, e^{\mu t}, e^{\nu t})$, 
where the values of $\lam$, $\mu$, $\nu$ are distinct, satisfy
\[ x''''=(\lambda+\mu+\nu)x'''-(\lambda\mu+\lambda\nu+\mu\nu)x''+
\lambda\mu\nu x', \]
where $\epsilon$ takes both $\pm 1$.
When, in addition, $\lambda+\mu+\nu=0$, we have $k=0$, $\epsilon=-1$,
$q=\sqrt{\lam^2+\lam\mu+\mu^2}$, and $M=\lam\mu(\lam+\mu)/q^3$.
}
\end{example}

\begin{example} \label{mk}
 {\rm For curves given by $x = (t, e^{\lam t}, te^{\lam t})$, the equation is
\[ x''''=2\lam x''' - \lam^2 x'',\]
and $\epsilon=-1$, $k=-\sqrt{2}{\rm sign}(\lam)$, $M=\sqrt{2}{\rm sign}(\lam)$.
In particular, the identity $M-k\epsilon=0$ holds.}
\end{example}

Together with these examples, the curves in the next table
exhaust the list of nondegenerate curves with constant general-affine
curvatures.  Note that  the
hyperbolic helix in Example \ref{helix} is listed also in the class
of curves numbered $1$, and Example \ref{logspiral} is in the class
numbered $7$ below. 
We list also the associated differential equations.
\bigskip

{
\centering
\begin{tabular}{lll} \hline
 & curves & differential equations \\ \hline
1& $(t,e^{\lambda t}, e^{\mu t})$ & $x''''=(\lambda +\mu)x'''-\lambda\mu x''$ \\
2& $(e^t, te^t, e^{\lam t})$ & $x''''=(\lam +2)x'''-(2\lam+1)x''+\lam x'$ \\
3& $(t,{1\over 2}t^2, e^{\lambda t})$ & $x''''=\lambda x'''$ \\
4& $(e^t, te^t, t^2e^t)$ & $x''''=3x'''-3x''+x'$ \\
5&$(t, e^t\cos(pt),e^t\sin(pt))$ & $x''''=2x'''-(p^2+1)x''$ \\
6&$(e^{\lam t},\cosh(pt), \sinh(pt))$ & $x''''=\lam x'''+p^2x''-\lam p^2 x'$ \\
7&$(e^{\lam t},e^{\mu t}\cos(pt),e^{\mu t}\sin(pt))$ &
  $x''''=(\lam + 2 {\mu})x'''-(p^2+{\mu (2\lam+\mu)})x''+\lam(p^2+{\mu}^2)x'$  \\
8 & $(t, {1\over 2}t^2, {1\over 6}t^3)$ & $x''''=0$ \\
\hline
\end{tabular}
}

\subsection{From equiaffine to general-affine for space curves} 
\label{3equiaffine}
Let us pay some attention to the equiaffine theory of space curves 
in comparison with the general-affine treatment.
Recall the choice of the frame $e=\{e_1,e_2,e_3\}$ and 
the scalars $\ell$ and $m$ in {Section} \ref{subsc:3frames}, 
\eqref{spccov}.
\begin{equation} \label{pfaffequiaffine}
 d \left(\begin{array}{c} x \\ e_1 \\ e_2 \\ e_3\end{array}\right)
=
\left(
\begin{array}{ccc}
\om & 0 & 0 \\
0   & \om & 0 \\
0   & 0   & \om \\
-m\om & -\ell\om & 0
\end{array}\right)
\left(\begin{array}{c} e_1 \\ e_2 \\e_3 \end{array}\right).
\end{equation}
In the equiaffine treatment, it is enough to consider only the
unimodular change of frames: {\it i.e.} $\lambda\nu\gamma=1$ 
and $e$ takes values in SL$(3,{\bf R})$.
By $(\ref{unimo})$, we have $\lambda^6=1$.
This means that the scalar $\ell$ is an absolute invariant
and the scalar $m$ is an invariant determined up to $\pm 1$
by $(\ref{spccov})$. As was remarked in Section \ref{subsc:3frames},
the scalar $\ell$ is usually 
called the \textsl{equiaffine curvature}
of the space curve and the scalar $m$ is 
called the \textsl{equiaffine torsion}; 
we refer to the books \cite{Bl,Sc}.
The invariant $\ell$ measures how the space curve differs from
the osculating cubic parabola, which is defined to be
the curve $(t,t^2/2, t^3/6)$ relative to certain affine coordinates.

The parameter $t$ for which $\om=dt$ holds 
is called an \textsl{equiaffine length parameter.}
Then 
the equation $(\ref{pfaffequiaffine})$ implies that
the immersion $x(t)$ satisfies the differential equation
\begin{equation} \label{easpcode}
x''''+ \ell x'' + mx' =0,
\end{equation}
which is written in the form of equation $(\ref{x3ode})$ where
$a=0$, $b=-\ell$, $c=-m$. By $(\ref{galength})$, the general-affine
length parameter $\sigma$ is determined by use of the equiaffine
curvature $\ell$ as 
\[ d\sigma^2 = Ldt^2,\qquad {\rm where}\quad L=\epsilon\ell,\quad
\epsilon={\rm sign}(\ell)\]
and, by $(\ref{kforspacecurve})$ and $(\ref{gasecond})$, 
the first and second general-affine curvatures are given as
\[
k=L'L^{-3/2}, \qquad M= mL^{-3/2}\]
in terms of equiaffine curvature and equiaffine torsion.
Relative to the parameter $\sigma$, the map $y(\sigma)=x(t)$ is seen to
satisfy the equation $(\ref{yode})$ whose coefficients are given by
\begin{eqnarray*}
A(\sigma) &=& -{3L'\over L^{3/2}}, \\
B(\sigma) &=& -\epsilon + {L'\over 4L^3} - {2L''\over L^2}, \\
C(\sigma) &=& L^{-3/2}\left(-m -{\epsilon L'\over 2}
- {L'''\over 2L}+{3L'L''\over 4L^2}-{3L'^3\over 8L^3}\right),
\end{eqnarray*}
by use of $(\ref{newspca})-(\ref{newspcc})$.

The curves for which $\ell$ and $m$ are constant 
can be listed, by solving the equation $(\ref{easpcode})$,
as follows; see \cite[p.75]{Sc}.  
\medskip

\begin{tabular}{lll}
1. $(e^{\lambda t}, e^{\mu t}, e^{-(\lambda+\mu)t})$, &
2. $(te^{\lambda t}, e^{\lambda t}, e^{-2\lambda t})$, &
3. $(e^{-2\alpha t}, e^{\alpha t}\cos(\beta t), e^{\alpha t}\sin(\beta t))$,\\
4. $(t, \cosh t, \sinh t)$, &
5. $(t, \cos t, \sin t)$, &
6. $(t,{1\over 2}t^2, {1\over 6}t^3)$, 
\end{tabular}

\medskip

\noindent where $\lambda$, $\mu$, $\alpha$, $\beta$ are nonzero constants.
They are homogeneous under equiaffine transformations.
The value $m$ is nonzero for the first three 
and is zero for the last three. The value $\ell$ is $-1$, $1$ and $0$ for
the last three, in this order.
Except for the last example, 
 the general-affine curvature $k$ is defined, and
it is vanishing because $\ell$ is constant.
The listed curves are general-affinely equivalent to some
of the examples in the previous section.

\subsubsection{Extremal equiaffine space curves}

W. Blaschke \cite{Bl} gave a variational formula of the equiaffine
length and showed that extremal curves of this variation are
the curves with $\ell=m=0$; hence, the cubic parabola.
This will be seen as follows in the present setting.

\begin{theorem}[\cite{Bl}]\label{equi3extremal}
A nondegenerate curve in the affine $3$-space  is extremal
relative to the equiaffine length functional if and only if
the equiaffine curvatures $\ell$ and $m$ are vanishing.
\end{theorem}

\noindent {Proof.} 
 Let $x_{\eta} (t)$ denote a family of curves parametrized 
by $\eta$ around $\eta=0$ and  $x_0 = x$ as before.
We assume that $x_{\eta} (t)= x(t)$ outside of a compact set $C$
 and $x_0(t)$ is parametrized by equiaffine arc length, and that
 $\omega$ does not vanish anywhere for all $\eta$. 
 The equiaffine length functional is given by  
\[
 L(\eta) = \int_C \omega(\eta).
\]
Then the curve $x$ is \textsl{equiaffine extremal}  if 
$ \delta L =0$.
 Let $e= \{x, e_1, e_2, e_3\}$ be the frame defined as in \eqref{pfaffequiaffine}, and set $\Omega$ to be $4 \times 3$ coefficient matrix.
 Since $\{e_1, e_2, e_3\}$ are linearly independent, there exists a 
 $4 \times 3$-matrix $\tau$ such that 
 \[
  \delta  
\begin{pmatrix}
 x \\ e_1 \\ e_2 \\  e_3
\end{pmatrix}
 = \tau 
\begin{pmatrix}
 e_1 \\ e_2 \\  e_3
\end{pmatrix}.
 \]
We denote the components of $\Omega$ and $\tau$ 
by $\omega_{\alpha}^{\beta}$ and $\tau_{\alpha}^{\beta}$, respectively,
where $0\le \alpha\le 3$ and  $1\le \beta\le 3$.
 Since $\delta d e = d \delta e$, 
 we have 
$
 \delta \omega_{\alpha}^{\beta} -d \tau_{\alpha}^\beta
 = \tau_{\alpha}^{\gamma} \omega_{\gamma}^{\beta}
- \omega_{\alpha}^{\gamma} \tau_{\gamma}^{\beta}
$; in terms of entries of $\Omega$ and $\tau$, we have
 \begin{align}
\delta \omega - d \tau_{0}^1  &=  - (m \tau_0^3 + \tau_1^1) \omega, 
 \label{eqaff0-1}\\
- d \tau_{0}^3  &=  (\tau_0^2 - \tau_1^3) \omega, \label{eqaff0-3}\\
\delta \omega - d \tau_{1}^2  &=  (\tau_1^1 - \ell \tau_1^3  - \tau_2^2) \omega, 
 \label{eqaff1-2}\\
\delta \omega - d \tau_{2}^3  &=  (\tau_2^2 - \tau_3^3) \omega. 
 \label{eqaff2-3}
\end{align}
 First, note that since $\{e_1, e_2, e_3\}$ takes values 
in ${\rm SL}(3,{\bf R})$, we have 
\begin{equation}\label{eq:tausum}
 \tau_1^1 + \tau_2^2 + \tau_3^3 =0.
\end{equation}
 Adding \eqref{eqaff0-1}, \eqref{eqaff1-2} and  \eqref{eqaff2-3}, we have
\[
 3 \delta \omega - d \tau_0^1 - d \tau_1^2 - d \tau_2^3 
= -(m \tau_0^3 + \ell \tau_1^3 + \tau_3^3)\omega.
\]
 On the one hand, subtracting \eqref{eqaff0-1} from  \eqref{eqaff2-3}, we have
\[
 - d \tau_2^3 + d \tau_0^1 = (m \tau_0^3 -2 \tau_3^3) \omega,
\]
 where we use the relation \eqref{eq:tausum}.
 Thus we have 
\[
 3 \delta \omega = d \tau_0^1 + d \tau_1^2 + d \tau_2^3 - (m \tau_0^3 + \ell \tau_1^3)
 \omega - \frac{1}{2} (d \tau_2^3 - d \tau_0^1 + m \tau_0^3 \omega),
\]
 and therefore 
\[
 3 \delta L = \int_C \left( - \frac{3}{2} m \tau_0^3 - \ell \tau_1^3\right) \omega
\]
 holds. Finally, using \eqref{eqaff0-3} and integration by parts, 
we obtain
\begin{align*}
 3 \delta L &= \int_C \left\{ - \frac{3}{2} m \tau_0^3 - \ell 
 \left({ \tau_0^3}' + \tau_0^2\right)
 \right\}\omega \\
 & = 
 \int_C \left\{ \left(- \frac{3}{2} m + \ell' \right) \tau_0^3 - \ell \tau_0^2
 \right\}\omega,
\end{align*}
 where the $\{'\}$ denotes the derivative with respect to the arc length. 
 Since $\tau_0^3$ and 
 $\tau_0^2$ are independent variation vector fields, we have 
 completed the proof.

\subsection{From general-affine to projective for space curves}
\label{remtoproj}

\comment{
Study of projective curves has an old history: a first systematic study
was done by Halphen \cite{Ha}, followingly by Laguerre and Forsyth, and then,
from the viewpoint of projective differential geometry, 
by Wilczynski \cite{Wi}.
Here, we present some invariants of projective space curves due
to \cite[p.83-84]{La} and show a relation of the projective length
parameter with the $\theta$ in \eqref{theta}; 
in Appendix, we will make some compliment on this matter.
}

A space curve in ${\bf P}^3$ is given by 
the immersion $t \longmapsto x(t)\in {\bf A}^4$ in homogeneous
coordinates satisfying an ordinary differential equation of the form
\[
x''''+4p_1x'''+6p_2x''+4p_3x'+p_4x=0.\]
By multiplying a nonzero factor to the indeterminate $x$, the equation is
transformed to the equation
\[
x'''''+6P_2x''+4P_3x'+P_4x=0,\]
where 
\[
\begin{array}{rcl}
P_2 &=& p_2 - p_1^2 -p_1',\\
P_3 &=& p_3 - 3p_1p_2+2p_1^3-p_1'', \\
P_4 &=& p_4-4p_1p_3+6p_1^2p_2-6p_1'p_2-3p_1^4+6p_1^2p_1'+3(p_1')^2-p_1'''.
\end{array}
\]
Then, the two forms $\theta_3dt^3$ and $\theta_4dt^4$, where
\begin{equation} \label{thetainv}
\begin{array} {rcl}
\theta_3 &=& P_3-{3\over 2}P_2',\\
\theta_4 &=& P_4 -{9\over 5}P_2''-{81\over 25}P_2^2-2\theta_3',
\end{array}
\end{equation}
are fundamental invariant forms: \cite{La}. 
 Provided that $\theta_3\neq0$, the parameter $s$ defined as
\[ ds = \theta_3^{1/3} dt\]
is called the \textsl{projective length parameter}. 
Relative to this parameter, 
we can define projective curvatures; we refer to Appendix \ref{subs-projinv}.
When $\theta_3\equiv 0$, the curve $x$ has a special property that
the curve formed by the tangent vectors to the curve $x$, which lies 
in the $5$-dimensional projective space
consisting of lines in ${\bf P}^3$, 
is degenerate in the sense that it belongs to a $4$-dimensional
hyperplane. Such a curve was said to belong to a \textsl{linear complex} and
is named Gewindekurve in \cite{Bl}. 

Given a nondegenerate curve $x(t)$ in the affine space ${\bf A}^3$, which is
described by the differential equation $(\ref{gaspcode})$, 
we associate a curve in ${\bf P}^3$ by a mapping $t \longmapsto 
(1,x(t))\in {\bf A}^4$, where $1$ is a constant function. 
Then, the projective invariants are computed by the definition above. 
In fact, a straightforward computation shows that
\begin{align}
& \theta_3 = {1\over 4}(M-\epsilon k), \label{theta3} \\
& \theta_4 = -{3\over 4}kM-{1\over 2}M'+ {1\over 5}\epsilon k'
 + {3\over 10}\epsilon k^2 -{9\over 100}. \label{theta4}
\end{align}
In particular, when $\theta_3\neq 0$,
the projective length parameter $s$ is given as above by use of the general-affine
curvatures $k$ and $M$.
When $M=\epsilon k$, the curve belongs to a linear complex.
Example \ref{mk} in the previous subsection
is such an example.

\section{General-affine extremal space curves and the associated
differential equations}
 In Section \ref{subsc:3frames}, we have defined a frame 
 for a general-affine space curve under the condition
that the curve has no affine inflection point.
In this section, we obtain the condition under which a space curve
is extremal relative to the length functional and, in particular,
show that any curve with constant general-affine curvatures is
extremal.

 Let $x_{\eta}(t)$ be a family of curves parametrized by 
 $\eta$ around $\eta = 0$ and $x_0 = x$. We assume that $x_{\eta} (t) = x(t)$ 
 outside a compact set $C$, and that the invariant $\omega_s$ does not vanish anywhere for all $\eta$. Then $x_{\eta}$ and the corresponding frame $\{E_1, E_2, E_3\}$
 satisfy the equation in \eqref{pfaff3}.
\comment{
\begin{equation} \label{pfaff3-2}
 d \left(\begin{array}{c} x_\eta \\ E_1 \\ E_2 \\ E_3\end{array}\right)
=\Omega
\left(\begin{array}{c} E_1 \\ E_2 \\E_3 \end{array}\right),
\qquad
\Omega =\left(
\begin{array}{ccc}
 \om_s    &   0  &  0 \\
-{1\over 2}k\om_s & \om_s & 0 \\
0 & -k\om_s  & \om_s \\
- M \om_s & -\epsilon \om_s & -{3\over 2}k\om_s
\end{array}\right).
\end{equation}
}
 Then the length functional $L$ is given by 
\[
 L(\eta) = \int_C \omega_s (\eta),
\]
 and the curve $x=x_0$ is said to be \textsl{general-affine  extremal} if 
\[
 \delta L = \frac{d L}{d \eta}\Big|_{\eta =0}=0
\]
 holds for any compactly supported deformation of $x$.


 We now consider the variation 
\[
 \delta  \begin{pmatrix}
x_{\eta}  \\ E_1 \\ E_2  \\ E_3
 \end{pmatrix} = \tau  \begin{pmatrix}
  E_1 \\ E_2  \\ E_3
\end{pmatrix}, 
 \quad \quad 
\tau = (\tau_\alpha^\beta)_{0 \leq \alpha\leq 3,  1\leq \beta \leq 3}.
\]
 Then the compatibility condition $d \delta = \delta d$ implies that
 \[
 \delta \omega_\alpha^\beta- d \tau_\alpha^\beta
 = \tau_\alpha^\gamma \omega_\gamma^\beta  - \omega_\alpha^\gamma \tau_\gamma^\beta,
 \]
 where we set the entries of $\Omega$ in \eqref{pfaff3} as  
 $(\omega_\alpha^\beta)_{0 \leq i\leq 3, 1\leq \beta \leq 3}$.
 Then they are explicitly given by 
\begin{align}
 \delta \omega_s - d \tau_0^1 &= \left(
 - \frac{1}{2} k \tau_0^1 - M \tau_0^3 - \tau_1^1 \right) \omega_s, \label{eq:01} \\
- d \tau_0^2   &=  (\tau_0^1 -k \tau_0^2 - \epsilon \tau_0^3 -\tau_1^2) \omega_s, \label{eq:02}\\
- d \tau_0^3   &=  \left(\tau_0^2 -\frac{3}{2}k \tau_0^3 -\tau_1^3 \right) \omega_s, \label{eq:03}\\
 \delta \omega_s - d \tau_1^2 &= \left(
 \tau_1^1 - \frac{1}{2}k \tau_1^2 - \epsilon \tau_1^3 - \tau_2^2\right) \omega_s, \label{eq:12}\\
- d\tau_1^3 &= \left(\tau_1^2 - k \tau_1^3 -\tau_2^3\right) \omega_s, \label{eq:13}\\
- d \tau_2^1 &= \left(\frac{1}{2} k \tau_2^1 - M \tau_2^3 - \tau_3^1 \right) \omega_s, \label{eq:21}\\
 \delta \omega_s - d \tau_2^3 &= \left(
\tau_2^2 -  \frac{1}{2} k \tau_2^3  - \tau_3^3 \right) \omega_s, \label{eq:23}\\
 - \epsilon \delta \omega_s - d \tau_3^2 &= \left(
 -\epsilon (\tau_3^3 - \tau_2^2 )+ \tau_3^1 + \frac{1}{2} k \tau_3^2 + M \tau_1^2
 \right) \omega_s, \label{eq:32}\\
- \frac{1}{2}\delta (k \omega_s) - d \tau_1^1  &= (- M \tau_1^3 -\tau_2^1) \omega_s,
\label{eq:11}\\
- \delta (k \omega_s) - d \tau_2^2  &= (\tau_2^1 - \epsilon \tau_2^3 - \tau_3^2) \omega_s,
 \label{eq:22}\\
- \frac{3}{2}\delta (k \omega_s) - d \tau_3^3  &= (\tau_3^2 + M \tau_1^3 + \epsilon 
 \tau_2^3) \omega_s, \label{eq:33} \\
- \delta (M \omega_s) - d \tau_3^1 &= (k \tau_3^1 + M (\tau_1^1 - \tau_3^3) + \epsilon 
 \tau_2^1) \omega_s. \label{eq:31}
\end{align}
 Adding \eqref{eq:23} and $-\epsilon$\eqref{eq:32}, we have 
\begin{equation*}
2 \delta \omega_s -  d \tau_2^3 + \epsilon d \tau_3^2 = 
\left(- \frac{1}{2} k (\tau_2^3  + \epsilon \tau_3^2) - \epsilon \tau_3^1 - \epsilon 
 M \tau_1^2\right) \omega_s.
\end{equation*}
Then, by Stokes' theorem, we have 
\begin{equation}\label{eq:omega}
 2 \delta \int_C \omega_s = 
 \int_C  \left(-\frac{1}{2}k(
 \tau_2^3 + \epsilon \tau_3^2) - \epsilon \tau_3^1 - \epsilon M \tau_1^2 \right)
 \omega_s. 
\end{equation}
 Next, from \eqref{eq:11} $+$ \eqref{eq:22} $-$\eqref{eq:33}, we get
\begin{equation*}
 (- {\tau_2^2} -{\tau_1^1} +{\tau_3^3})' = -2 (\epsilon \tau_2^3 + \tau_3^2)
 - 2 M \tau_1^3,
\end{equation*}
 which is written as
\begin{equation}\label{eq:tau2332}
 -\frac{1}{2}k(
 \tau_3^2 + \epsilon \tau_2^3) = 
 - \frac{1 }{4}\epsilon k ({\tau_1^1} + {\tau_2^2} - {\tau_3^3})' + \frac{1}{2} \epsilon k M \tau_1^3.
\end{equation}
 Here the $\{'\}$ denotes $\frac{d}{\omega_s}$.
 On the one hand, by \eqref{eq:21}, 
\begin{equation}\label{eq:21-2}
  \tau_3^1 = {\tau_2^1}' + \frac{1}{2} k \tau_2^1 - M \tau_2^3
\end{equation}
 and \eqref{eq:22}$+$\eqref{eq:33}$-5$\eqref{eq:11} implies that
\begin{equation}\label{eq:223311}
  (-{\tau_2^2} -{\tau_3^3} + 5{\tau_1^1})' = 6 (\tau_2^1  + M \tau_1^3).
\end{equation}
 Then by use of \eqref{eq:13} and \eqref{eq:223311}, 
\eqref{eq:21-2} can be  rephrased as 
\begin{equation}\label{eq:tau31}
 \tau_3^1 = {\tau_2^1}' + \frac{1}{12} k ( 5 \tau_1^1 - \tau_2^2 - 
 \tau_3^3)' - M\left({\tau_1^3}' + \tau_1^2 -\frac{1}{2}k \tau_1^3\right).
\end{equation}
 Finally \eqref{eq:tau2332} and \eqref{eq:tau31} implies that 
\begin{align}
 2 \delta \int_C \omega_s &= 
 \frac{\epsilon}{12} \int_C  
 \left(- k \left(8{\tau_1^1} + 2 \tau_2^2
 - 4 \tau_3^3\right)' + 12 M {\tau_1^3}'\right) \omega_s \nonumber \\
 & = \frac{\epsilon}{12} \int_C  
 \left(k' \left(8\tau_1^1 + 2 \tau_2^2 - 4 \tau_3^3\right) 
 + 12 M {\tau_1^3}'\right) \omega_s. 
 \label{eq:domegas}
\end{align}
 Here we use integration by parts for the second equality.

 We now compute $-6$\eqref{eq:01}$+2$\eqref{eq:12}$+4$\eqref{eq:23}. 
 A straightforward computation shows that 
\begin{align}
8 \tau_1^1 + 2 \tau_2^2 -4 \tau_3^3
 &= (6 \tau_0^1- 2\tau_1^2 -4\tau_2^3)'
 -3 k \tau_0^1 -6 M \tau_0^3 + k \tau_1^2 +2 \epsilon \tau_1^3 + 2 k \tau_2^3 
 \nonumber \\
 & = 6 X' -4 Y'
 -3 k X + 2 k Y -6 M \tau_0^3 + 2 \epsilon \tau_1^3. \nonumber
\end{align}
 Here $X = \tau_0^1 - \tau_1^2$ and $Y =\tau_2^3 - \tau_1^2$.
 Thus \eqref{eq:domegas} can be again rephrased, by using integration by 
 parts, as 
\begin{align}
 24 \epsilon \delta \int_C \omega_s =
\int_C 
\left\{ (-6 k''-3 k' k) X + (4 k'' + 2 k' k) Y 
\right. & -6 k' M \tau_0^3 \nonumber \\
 &\left.+ (2 \epsilon k'- 12 M') \tau_1^3\right\} \omega_s.
\label{eq:domegas-f}
\end{align}
 Then by \eqref{eq:13} and \eqref{eq:02} we have 
 \[
  X = \tau_0^1 - \tau_1^2 = -{\tau_0^2}' + k \tau_0^2 + \epsilon \tau_0^3, \quad 
  Y = \tau_2^3 - \tau_1^2 = {\tau_1^3}' - k \tau_1^3.
 \]
 Finally, making use of \eqref{eq:03} to erase the $\tau_1^3$-term,
 we can see that the $\tau_0^2$ part of the integrand of \eqref{eq:domegas-f}
 is computed as 
\begin{equation}\label{eq:tau02part}
 -10 k''' -15 k'' k -5 k' k^2 -5 {k'}^2 +2 \epsilon k' -12 M'.
\end{equation}
 Similarly the $\tau_0^3$ part of the integrand of \eqref{eq:domegas-f}
 can be computed as
\begin{equation}\label{eq:tau03part}
 4 k'''' + 12 k''' k + (11 k^2 + 10 k' -8 \epsilon ) k'' + 
 7 {k'}^2 k -6 \epsilon k' k + 3 k' k^3 -6 k' M + 12 M'' + 18 
 M' k.
\end{equation}    
\begin{theorem} \label{spaceextremal}
A {nondegenerate} space curve without affine inflection point
is general-affine extremal
if and only if the following pair of ordinary differential equations 
 is satisfied$:$ 
\begin{align}
&k''' + \frac{3}{2}k k'' +\frac{1}{2} {k'}^2 + \frac{1}{2}k^2 k'  
 - \frac{1}{5} \epsilon k' + \frac{6}{5} M'=0 \label{eq:tau02part2}\\
\intertext{and}
&{k}'' + \frac{2}{3} k' k + \frac{5 }{6}\epsilon  k' M 
 -\frac{3 }{2} \epsilon k M' - \epsilon {M}''=0.
\label{eq:tau03part2}
\end{align}
In particular, all space curves which have constant 
general-affine curvatures are general-affine extremal.
\end{theorem}
\noindent Proof.
 Inserting \eqref{eq:tau02part}$=0$ into \eqref{eq:tau03part}$=0$, we have the 
 differential equation \eqref{eq:tau03part2}.

\begin{example} \label{mconst}
Extremal curves with constant $M$.
{\rm
First, assume $M=0$. Then 
\eqref{eq:tau03part2} can be easily integrated as
\[ k(t)= -3a\tan(at)\quad {\rm and}\quad 3a\tanh(at),\]
where $a$ is constant. Inserting this expression into 
\eqref{eq:tau02part2}, we get solutions 
\[ k(t)= -3a\tan(at),\quad a=\sqrt{2/5} \quad {\rm when} \quad \epsilon=1\]
and
\[ k(t)= 3a\tanh(at),\quad a=\sqrt{2/5} \quad {\rm when} \quad \epsilon=-1.\]
Second, assume $M$ is a nonzero constant. 
 Then
\[k(t)=-\frac{5}{4}\epsilon M + 3a \tanh(at)\]
is a solution of \eqref{eq:tau03part2} and it satisfies
\eqref{eq:tau02part2} if and only if
\[ a^2(80 a^2 - 125 \epsilon^2 M^2 + 32\epsilon)=0.\]
Thus, except for a constant solution, we have 
the above $k(t)$, where $a=\sqrt{-32 + 125 M^2}/(4 \sqrt{5})$ when $\epsilon=1$ and
$a=\sqrt{32+125M^2}/(4 \sqrt{5})$ when $\epsilon=-1$.
If we started with $-5/4\epsilon M-3a\tan(at)$, an another solution of
\eqref{eq:tau03part2}, then $a$ turns out to be pure imaginary and
we get the same curvature function.
}
\end{example}

We here recall the invariant $\theta_3$ given by the equation
$(\ref{theta3})$:
\begin{equation*} 
\theta_3 = {1\over 4}(M - \epsilon k).
\end{equation*}
Then, the differential equations $(\ref{eq:tau02part2})$ and
$(\ref{eq:tau03part2})$ are written as
\begin{align}
&k'''+{3\over 2}kk''+{1\over 2}k'^2+{1\over 2}k^2k'+\epsilon k'
 + {24\over 5}\theta_3' =0,  \\
\noalign{\smallskip}
& \theta_3'' + {3\over 2}k\theta_3' - {5\over 6}k'\theta_3 =0.
\end{align}
Since $\theta_3=0$ characterizes a curve belonging to a linear complex, 
see Section \ref{remtoproj}, in view of the equation 
\eqref{eq:variationformula}, we have the following corollary.

\begin{corollary} \label{lincomplex}
 The general-affine extremal space curve belongs to a linear complex 
 if and only if $M = \epsilon k$ and 
 $k$ satisfies \eqref{eq:variationformula}.
\end{corollary}

Since the differential equation 
\eqref{eq:variationformula} is the equation for the plane extremality,
we have the following method of constructing an extremal space
curve belonging to a linear complex:

\begin{corollary} \label{planetospace} 
Let $k$ be the general-affine curvature of an extremal plane curve without
affine inflection point. Let $\epsilon$ denote the
sign of this curve. Then, the set $\{k, M, \epsilon\}$, where $M=\epsilon k$,
defines a space curve that is general-affine extremal 
 and belonging to a
linear complex.
\end{corollary}

\def\erf{{\rm erf}}
Thanks to Example \ref{planeextremalexample}, 
we can give concrete examples of such curves in
Corollary \ref{planetospace}.
The explicit integration of the associated differential equation
can be carried out with computer assistance. For example,
when $\epsilon=-1$ and $k(t)=\sqrt{2}+3/t$, we get the curve
$(x_1, x_2, x_3)$ for $t>0$, where
\begin{align*}
x_1 &= {1\over t}, \\
x_2 &= 2^{1/4}\sqrt{\pi}\ \erf\left({\sqrt{3t} \over2^{1/4} }\right)
   - {1- \sqrt{2}\ 3t \over (3t)^{3/2}} \exp\left(- {3\over \sqrt{2}} t\right), \\
x_3 &= \int \left\{\frac{6}{t^{2}} \int H(t) \ dt
    +{1 \over t^{5/2}(\sqrt{2}+6 t)} \exp\left(- {3\over \sqrt{2}}t\right)
\right\}  \ dt, \\
\end{align*}
 with $\erf(x) = \displaystyle {2 \over \sqrt{\pi}} \int_0^{x} e^{-t^2} dt$ and
$H(t)=\displaystyle  \frac{1}{\sqrt{t}(\sqrt{2}+6t)^{2}}
 \exp\left(-{3\over \sqrt{2}} t\right) $. 

\newpage


\appendix
\noindent \begin{center} {\large \bf Appendix} \end{center}
\section{Projective treatment of plane curves}  \label{projtreat}

A study of basic notions such as
projective length parameter and projective curvature of
projective plane curves by use of moving frames was originally 
done by  E. Cartan in \cite{Ca1}, \cite[Chapitre 2]{Ca2};
nonetheless, in this appendix, we recall how to define 
such notions in a fairly detailed way 
because to the authors' knowledge it may not be so familiar in the present day
and it is useful for comparison with the general-affine treatment.

By a  projective plane curve we mean a nondegenerate immersion into 
${\bf P}^2$: $t\longrightarrow \underline{x}(t)\in {\bf P}^2$. 
We denote its lift to the affine space ${\bf A}^3$ by 
$t\longrightarrow x(t)\in {\bf A}^3-\{0\}$.

Let $e=\{e_0, e_1,e_{2}\}$ be a frame along $x$; at each 
point of $x$ it is a set of independent vectors of ${\bf A}^{3}$
which depends smoothly on the parameter. We choose $e_0=x$ for simplicity.
Then, the vector-valued $1$-form $de_0$ is written as
\[de_0 = \om_0^0e_0+ \om_0^1 e_1 + \om_0^2 e_2,\]
and the dependence of $e_i$ on the parameter is similarly written as 
\[ de_i = \sum_{j=0}^2 \om_i^j e_j,\qquad i=1,2.\]
In the following, we consider frames
such that the space generated by $e_0$ and $e_1$ is the space generated
by the vector $x(t)$ and the tangent vector $x'(t)$.
Then, $\om:=\om_0^1$ is nontrivial and $\om_0^2=0$. 
We consider furthermore frames such that the
condition $\om_0^0+\om_1^1+\om_2^2=0$ always holds.
Thus, we have the equation
\begin{equation} \label{streq}
\left(\begin{array}{c} d{e}_0 \\ d{e}_1 \\ d{e}_2 \end{array}\right) 
= \left(\begin{array}{ccc} \om_0^0 & \om & 0 \\
\om_1^0 & \om_1^1 & \om_1^2 \\ 
\om_2^0 & \om_2^1 & \om_2^2 \end{array}\right)
\left(\begin{array}{c} {e}_0 \\ {e}_1 \\ {e}_2 \end{array}\right),
\end{equation}
where the matrix $(\om_i^j)$ is called the coframe. 

\comment{To speak more abstractly, any frame is a smooth section of
the induced bundle along $\underline{x}$ of the frame bundle
of ${\bf P}^2$ with PGL$(2)$ as the group of transformations
and $(\om_i^j)$ is the induced Maurer-Cartan form on the bundle.}

Once we choose a frame $e$ with the required property, 
then another frame $\tilde{e}$ is given as
\begin{equation} \label{pframe}
\left(\begin{array}{c} \tilde{e}_0 \\ \tilde{e}_1 \\ \tilde{e}_2 \end{array}
\right) = \left(\begin{array}{ccc} \lambda & 0 & 0 \\
\mu & \alpha & 0 \\ \nu & \beta & \gamma \end{array}\right)
\left(\begin{array}{c} {e}_0 \\ {e}_1 \\ {e}_2 \end{array}\right),
\end{equation}
where $\lambda\alpha\gamma=1$.

We write the coframe for the frame $\tilde{e}$ by $\tilde{\om}_i^j$.
Since the derivation $d\tilde{e}_0$ has two expressions
$d\te_0 = \tom_0^0\te_0 + \tom\te_1$ and
$d\te_0 = d(\lam e_0) = d\lam e_0 + \lam(\om_0^0 e_0 + \om e_1)$,
we get the identities
\begin{eqnarray}
&\alpha \tom = \lam \om, & \label{ch1} \\
& \lam \tom_0^0 + \mu \tom = \lam \om_0^0 + d\lam .\label{ch2}&
\end{eqnarray}
Similarly considering the frames $\te_1$ and $\te_2$, we get
\begin{eqnarray}
& \gamma\tom_1^2 = \alpha \om_1^2, & \label{ch3} \\
& \alpha\tom_1^1+\beta\tom_1^2 = \alpha\om_1^1 +\mu\om + d\alpha, 
& \label{ch4}\\
& \lam \tom_1^0 + \mu\tom_1^1+\nu\tom_1^2 
  = \alpha\om_1^0 + \mu\om_0^0 + d\mu, & \label{ch5} \\
&\gamma\tom_2^2 = \gamma\om_2^2 + \beta\om_1^2+d\gamma, & \label{ch6} \\
&\alpha\tom_2^1+\beta\tom_2^2 = \gamma\om_2^1+\beta\om_1^1
+\nu\om + d\beta,& \label{ch7} \\
&\lam\tom_2^0+\mu\tom_2^1+\nu\tom_2^2
 = \gamma\om_2^0+\beta\om_1^0+\nu\om_0^0+d\nu. & \label{ch8}
\end{eqnarray}

\subsection{Projective curvature of plane curves}\label{subsc:Pframe}
Making use of these identities, we reduce the freedom of choice
of frames by a stepwise procedure. First, by \eqref{ch1},
\[ \tom = {\lam\over \alpha}\om.\]
If we set $\om_1^2=h\om$ and $\tom_1^2=\tilde{h}\tom$, then 
by \eqref{ch3}, we have
\[\tilde{h}=\alpha^3 h.\]
We assume that the curve is nondegenerate: $h\neq 0$. 
Then, by a choice of $\alpha$, we see that
there is a frame with $\tilde{h}=1$. 
Then, we can restrict our consideration to frames with $h=1$
and, therefore we necessarily have that $\alpha=1$ and $\lam\gamma=1$. 
This means that 
\[ \om_1^2=\om,\qquad \tom=\lam\om, \qquad \tom_1^2=\lam \om.\]
By \eqref{ch2}, we can choose $\mu$ so that $\tom_0^0=0$ and
to keep this condition for all frames, $\om_0^0=0$, we have
\[d\lam = \lam\mu \om.\]
Also by \eqref{ch6}, we can assume $\om_2^2=0$ and
\[d\gamma = -\beta\om.\]
Since $\om_0^0+\om_1^1+\om_2^2=0$ by assumption, we conclude
that $\om_1^1=0$, and by \eqref{ch4} we see that 
\[\mu=\lam \beta.\]
Now \eqref{ch5} reduces to $\tom_1^0=\gamma\om_1^0-\nu\om+\gamma d\mu$
and \eqref{ch7} reduces to $\tom_2^1=\gamma\om_2^1+\nu\om+d\beta$;
therefore
\[\tom_1^0-\tom_2^1 = \gamma(\om_1^0-\om_2^1)-2\nu\om + \gamma d\mu -d\beta.\]
Then, by a choice of $\nu$, there exists a frame with $\tom_1^0=\tom_2^1$.
Assuming this identity for all frames, we necessarily have
\[2\nu\om = \gamma d\mu - d\beta.\]
Using the formulas of $d\lam$ and $d\gamma$ above, 
and the identity $\lam\mu=1$, we see that
\[2\nu = \mu\beta.\]

Finally, by simplifying \eqref{ch8} using \eqref{ch5}, we see that
\[\lam \tom_2^0 = \gamma \om_2^0.\]
We set 
\[\om_2^0=\rho\om\]
and $\tom_2^0=\tilde{\rho}\tom$, and then
we have $\tilde{\rho}=\gamma^3\rho$.
Here we have two cases $\rho=0$ and $\rho\neq 0$.
In the latter case, we can find a frame with $\rho=-1$ and, therefore, we can
restrict the change of frame to the case $\gamma=1$. In this case,
by the formulas already obtained above, we see that
$\lam=1$, $\mu=\nu=\beta=0$; namely, we have uniquely determined 
the form $\om$ and the identity \eqref{ch5}
implies $\tom_1^0 = \om_1^0$. If we set
\begin{equation} \label{pcurv}
\om_1^0 = -k_p\om,
\end{equation}
and $\tom_1^0=-\tilde{k}_p\tom$, then it turns out that $\tilde{k}_p=k_p$.

Now, we call $\om$ the \textsl{projective length element} and $k_p$ the
\textsl{projective curvature}.
Thus we have the following:

\begin{lemma} \label{plemma}
Assume that the curve is nondegenerate and the
scalar $\rho$ is nonvanishing. Then, the frame is uniquely defined 
such that the coframe has the form
\begin{equation} \label{proplanepfaff}
\Omega =
\left(\begin{array}{ccc} 0 & \om & 0 \\
-k_p\om & 0 & \om \\ 
-\om & -k_p\om &  0 
\end{array}\right),
\end{equation}
where $k_p$ is the projective curvature and the $1$-form $\om$ is
the projective length element.
\end{lemma}
E. Cartan \cite{Ca2} called this coframe the formula of Frenet.

When $\rho=0$, the form $\om_1^0$ is
not uniquely determined and we need a separate consideration.
The condition $\rho=0$ is equivalent to $P=0$, 
 where $P$ is defined
in \eqref{projlen} in Section \ref{gatoproj}.

\subsection{Extremal projective plane curves}

 In this subsection, we derive the differential equation of 
 an extremal projective plane curve according to \cite{Ca1}.
 Let $x_{\eta} (t)$ be a family of projective plane curves 
 such that $x_0 = x$, and let $\omega$ be the projective length element 
 of $x_t$. Assume that $x= x_0$ is parametrized by projective arc 
 length, $x_{\eta} (t) = x(t)$ outside a compact set $C$ and 
 $\omega$ does not vanish everywhere for all $\eta$. 
 Consider the length functional 
\[
 L = \int_C \omega(\eta).
\]
 It is natural to call the curve $x=x_0$ \textsl{projective extremal} with 
 respect to the length functional if 
$ \delta L =0 $
 holds. In the following we compute the associated differential equation 
 for an extremal projective plane curve. 
 Let $k$ denote the curvature $k_p$ in this subsection.
 Let $e=\{e_0, e_1, e_2\}$ be the frame for a family of projective plane curves 
 defined as in \eqref{proplanepfaff}.
\comment{
in Section \ref{subsc:Pframe}, {\it i.e.}
\[
d  \begin{pmatrix}e_0 \\e_1 \\e_2    \end{pmatrix} 
 = \Omega\begin{pmatrix}
 e_0 \\e_1 \\e_2   
\end{pmatrix}.
\]
}
 Then the variation of the frame $e$ can be defined as 
\[
 \delta \begin{pmatrix}e_0 \\e_1 \\e_2    \end{pmatrix}  = \tau \begin{pmatrix}e_0 \\e_1 \\e_2    \end{pmatrix}, 
\]
 where $\tau$ is a $3\times 3$ matrix. Denote the entries of $\Omega$ and $\tau$ by 
 $\omega_{\alpha}^{\beta}$ and $\tau_{\alpha}^{\beta}$. Here the Greek letters run
 from $0$ through $2$. Then the compatibility $d \delta e = \delta d e$
 implies that 
\begin{align}
 \delta \omega - d \tau_0^1 &= (\tau_0^0- k \tau_0^2 - \tau_1^1) \omega,\label{eqp:01} \\
 - d \tau_0^2 &= (\tau_0^1- \tau_1^2) \omega,\label{eqp:02} \\
  - d \tau_1^1 &= (\tau_1^0- k \tau_1^2 + k \tau_0^1 - \tau_2^1)\omega, \label{eqp:11} \\
 \delta \omega - d \tau_1^2 &= (\tau_1^1+ k \tau_0^2 - \tau_2^2) \omega,\label{eqp:12} \\
 -\delta \omega - d \tau_2^0 &= (-k \tau_2^1- \tau_2^2 + \tau_0^0 + k \tau_1^0) \omega.\label{eqp:20} 
\end{align}
 We first note that since the frame $e=\{e_0, e_1, e_2\}$ takes values in 
 ${\rm SL}(3,{\bf R})$,
\begin{equation}\label{eq:sumtau}
 \tau_0^0 + \tau_1^1 + \tau_2^2 =0
\end{equation}
 holds. Now adding \eqref{eqp:01}, \eqref{eqp:12} and $-1$\eqref{eqp:20}, 
 we have 
\begin{equation}\label{eqp:3delta}
3 \delta \omega - d \tau_0^1 - d \tau_1^2 + d \tau_2^0 = k (\tau_2^1 - \tau_1^0)\omega.
\end{equation}
Then by \eqref{eqp:11} and \eqref{eqp:02}
\begin{equation}\label{eqp:tau21}
 \tau_2^1 - \tau_1^0 = {\tau_1^1}' -k {\tau_0^2}'
\end{equation}
 holds. Here the $\{'\}$ denotes the derivation 
with respect to the projective 
 length element, {\it i.e.} $a' = \frac{d a}{\omega}$ for a function $a$.
 Moreover, subtracting \eqref{eqp:01} from \eqref{eqp:12} and  
 using \eqref{eq:sumtau}, we have 
\[
 d (\tau_0^1 - \tau_1^2) = (3 \tau_1^1 + 2 k \tau_0^2) \omega.
\]
 Rephrasing the above equation by using \eqref{eqp:02}, we have 
\begin{equation}\label{eq:tau11}
 3 \tau_1^1 = -2 k \tau_0^2 + ({\tau_0^1}- {\tau_1^2})' 
 = -2 k \tau_0^2 - {\tau_0^2}''. 
\end{equation} 
 Inserting \eqref{eq:tau11} and \eqref{eqp:tau21} into \eqref{eqp:3delta},
 we have 
\begin{equation}\label{eqp:deltaomega}
 3 \delta \omega = d \tau_0^1 + d \tau_1^2 - d \tau_2^0 
 +\left(-{2\over 3}k{k}'\tau_0^2 - {5\over 3}k^2{\tau_0^2}'
-{1\over 3}k {\tau_0^2}'''\right)\om.
\end{equation}
 Finally, by integration by parts and Stokes' theorem, 
\[
 9 \delta L 
 =  \int_C \left(k''' +8 k k'\right) \tau_0^2 \omega.
\]
 Therefore we have the following:
\begin{theorem} \label{cartanplane}
{\rm \cite{Ca1}}
 A plane curve without inflection points is projective extremal 
 relative to the length functional 
 if and only if 
\[
 k''' +8 k k' =0
\]
 holds.
\end{theorem}

\section{Projective treatment of space curves} \label{appproj3}

We give a summary on how to define two kinds of projective curvatures of
space curves in projective $3$-space. We recall the normalization
of frame by use of the Halphen canonical form
of the differential equation, \cite{La},
and another normalization by G. Bol \cite{Bol}.
Then, following M. Kimpara \cite{Ki}, we compute the 
variational formula of the projective length functional
for the curve with $\theta_3\neq 0$, according to 
which we can see that such a curve is extremal
if and only if both projective curvatures are constant;
such curves are classified in Appendix C.

\subsection{Projective curvatures for space curves}  \label{subs-projinv}

In Section \ref{remtoproj}, we introduced the differential equation
\begin{equation} \label{Pdiff}
x'''''+6P_2x''+4P_3x'+P_4x=0,
\end{equation}
which describes any nondegenerate space curve in projective space
: $t\longrightarrow x(t)\in {\bf P}^3$, and 
defined two invariant forms $\theta_3dt^3$ and $\theta_4dt^4$ 
in \eqref{thetainv}.
We now assume $\theta_3\neq 0$ and choose the parameter $t$ so that
$\theta_3=1$; namely, let $t$ be a 
projective length parameter. Then, the equation is written as
\begin{equation} \label{Preddiff}
x''''+6P_2x'' + 2(2+3P_2')x'+P_4x=0,
\end{equation}
which is called the \textsl{Halphen canonical form}, and
the two scalars $P_2$ and $P_4$ (or $\theta_4$ instead of $P_4$) 
are called the \textsl{projective curvatures}; we refer to \cite[p.26]{FC}. 
For the sake of later reference, we set
\begin{eqnarray}
k &=& {3\over 5}P_2, \label{spcurvature}\\
\theta &=& \theta_4, \label{sptheta}
\end{eqnarray}
and call them the \textsl{first projective curvature} 
and the \textsl{second projective curvature}, respectively.
The curve with constant curvatures is called an \textsl{anharmonic}
curve and we give a complement to the study 
in \cite[Section 3 of Chapter 14]{Wi} by giving
a classification of such curves in Appendix \ref{appclass}.

With this preparation, let us choose a frame
$(e_1,e_2,e_3,e_4)=(x,x',x'',x''')$; then the coframe is given as
\begin{equation} \label{pfaffproj}
 d \left(\begin{array}{c} e_1 \\ e_2 \\ e_3 \\e_4 \end{array}\right)
=
\left(
\begin{array}{cccc}
0    &   1  &  0 & 0\\
0    &   0   & 1 & 0 \\
0    &   0   &  0 & 1 \\
-P_4 & -4-6P_2' & -6P_2 & 0 
\end{array}\right) dt 
\left(\begin{array}{c} e_1 \\ e_2 \\ e_3 \\e_4 \end{array}\right).
\end{equation}

On the other hand, in the book \cite[Section 39]{Bol}, another
choice of frame was given by directly using the invariants $\theta_3$
and $\theta_4$. 
It is done, for the differential equation \eqref{Pdiff},
by choosing $\{x, u, y, z\}$, where
\[u=x',\quad y=u'+{9\over 5}P_2 x,\quad
z=y'+{12\over 5}P_2u + 2\theta_3 x.\]
In terms of the frame $\{e_1,e_2,e_3,e_4\}$, it is given as
\begin{equation}\label{eq:pscf}x=e_1,\ u=e_2,\ y=e_3+{9\over 5}P_2e_1,
z=e_4+{21\over 5}P_2e_2+\left({9\over 5}P_2'+2\theta_3\right)e_1.
\end{equation}
Then, we can see that
\[
d\left(
\begin{array}{c} x\\ u\\y\\z\end{array}\right) =
\left(
\begin{array}{cccc}
0 & 1 & 0 & 0\\ -3k & 0 & 1 & 0\\
-2\theta_3 & -4k & 0 & 1 \\ -\theta_4 & -2\theta_3 & -3k & 0
\end{array}\right) dt 
\left(
\begin{array}{c} x\\ u\\y\\z\end{array}\right).
\]

When this choice of the frame 
is applied to the equation \eqref{Preddiff}, the frame equation is
\begin{equation} \label{bolframe}
d\left(
\begin{array}{c} x\\ u\\y\\z\end{array}\right) =
\left(
\begin{array}{cccc}
0 & 1 & 0 & 0\\ -3k & 0 & 1 & 0\\
-2 & -4k & 0 & 1 \\ -\theta & -2 & -3k & 0
\end{array}\right) dt 
\left(
\begin{array}{c} x\\ u\\y\\z\end{array}\right).
\end{equation}

\subsection{Extremal projective space curves}
We now derive the differential equation of 
 an extremal projective space curve according to 
 \cite{Ki}. Let $x_{\eta} (t)$ be a family of 
 projective space curves such that $x_0 = x$,
 and let $\omega$ be the projective length element of $x_t$.
 Assume that $x = x_0$ is parametrized by projective arc 
 length, $x_{\eta}(t) = x(t)$ outside a compact set $C$,   
 and that $\omega$ does not vanish anywhere for all $\eta$.
 Let $k_1$ and $k_2$, respectively, denote the curvatures 
 $k$ and $\theta$ 
 in this subsection. Let $e = \{e_1, e_2, e_3, e_4\}$
 be the frame for a family of projective space curves defined in 
 \eqref{eq:pscf}. Then, according to \eqref{bolframe}, we have 
\[
de = \left(
\begin{array}{cccc}
0 & \om & 0 & 0 \\ -3 k_1 \om & 0 & \om & 0 \\ 
-2 \om & - 4 k_1\om & 0 & \om \\
- k_2\om & -2 \om & -3 k_1\om & 0
\end{array}
\right) e.
\]
 The variation of $e$ can be computed as $\delta e = (\tau_i^j) e$.
 Then the compatibility condition $d \delta e = \delta 
 d e $ is equivalent to 
\begin{align}
 - d{\tau_0^0} &= (-3 k_1 \tau_0^1- 2 \tau_0^2 - k_2 \tau_0^3 - \tau_1^0)\omega \label{eq:ps00}\\ 
 \delta \om - d{\tau_0^1}&=(\tau_0^0 -4 k_1 \tau_0^2 -2 \tau_0^3 - \tau_1^1)\omega \label{eq:ps01}\\
 - d {\tau_0^2}&= (\tau_0^1 - 3 k_1 \tau_0^3 - \tau_1^2) \omega\label{eq:ps02}\\
 -d{\tau_0^3}&= (\tau_0^2 - \tau_1^3)\omega \label{eq:ps03}\\
-3\delta ( k_1 \om) - d{\tau_1^0} &=(3 k_1 \tau_0^0 - 3 k_1 \tau_1^1 -2 \tau_1^2
 -k_2 \tau_1^3 - \tau_2^0)\omega \label{eq:ps10}\\
-d {\tau_1^1} &= (3 k_1 \tau_0^1 + \tau_1^0 -4k_1 \tau_1^2 -2 \tau_1^3 - \tau_2^1)\omega \label{eq:ps11}\\
 \delta \om - d {\tau_1^2}&= (3 k_1 \tau_0^2 + \tau_1^1 - 3 k_1 \tau_1^3 - \tau_2^2) \omega\label{eq:ps12}\\
- d {\tau_1^3} &= (3 k_1 \tau_0^3 + \tau_1^2 - \tau_2^3) \omega\label{eq:ps13}\\
-2 \delta \om - d{\tau_2^0} &= (2\tau_0^0 +4 k_1 \tau_1^0 -3 k_1 \tau_2^1
 -2 \tau_2^2 - k_2 \tau_2^3 - \tau_3^0) \omega\label{eq:ps20}\\
- 4\delta ( k_1 \om)- d{\tau_2^1} &= (2\tau_0^1 +4 k_1 \tau_1^1 + \tau_2^0 -4 k_1 \tau_2^2 -2 \tau_2^3- \tau_3^1) \omega \label{eq:ps21}\\
- d{\tau_2^2} &= (2\tau_0^2 +4 k_1 \tau_1^2 + \tau_2^1 - 3 k_1 \tau_2^3 - \tau_3^2) \omega\label{eq:ps22}\\
\delta \om -d {\tau_2^3} &= (\tau_0^0+2 \tau_0^3 + \tau_1^1 +4 k_1 \tau_1^3 + 2 \tau_2^2) \omega \label{eq:ps23}\\
-\delta (k_2 \om) - d {\tau_3^0} &= (2 k_2 \tau_0^0 +2 \tau_1^0+k_2 \tau_1^1 
 +3 k_1 \tau_2^0 +k_2 \tau_2^2 - 3 k_1 \tau_3^1 -2 \tau_3^2) \omega \label{eq:ps30}\\
-2\delta \om - d {\tau_3^1} &= (2 \tau_0^0 + k_2 \tau_0^1 +4 \tau_1^1+ 3k_1 \tau_2^1  +2 \tau_2^2  + \tau_3^0 -4 k_1 \tau_3^2) \omega \label{eq:ps31}\\
-3\delta ( k_1 \om) - d {\tau_3^2} &= (3 k_1 \tau_0^0 +k_2\tau_0^2 + 3 k_1 \tau_1^1 +2 \tau_1^2 +6 k_1 \tau_2^2 + \tau_3^1)  \omega \label{eq:ps32}\\
- d {\tau_3^3}&= (k_2 \tau_0^3 +2 \tau_1^3+ 3 k_1\tau_2^3 + \tau_3^2) \omega. 
 \label{eq:ps33}
\end{align}
 Here we use the relation $\tau_3^3 = - \tau_0^0 - \tau_1^1 - \tau_2^2$.
First adding \eqref{eq:ps02} and \eqref{eq:ps13},
\begin{equation*}
\tau_0^1 - \tau_2^3 = - {\tau_0^2}' - {\tau_1^3}',
\end{equation*}
 and by using \eqref{eq:ps03}, the above equation can be rephrased as
\begin{equation}\label{eq:ps0123}
\tau_0^1 - \tau_2^3 = - 2{\tau_0^2}' - {\tau_0^3}''.
\end{equation}
 Here the $\{'\}$ denotes the derivation with respect to the arc length.
 Next by \eqref{eq:ps11}
\begin{align*}
\tau_1^0 - \tau_2^1 &= - {\tau_1^1}' +2\tau_1^3- k_1 (3  \tau_0^1 - 4 \tau_1^2) \\
 & =  - {\tau_1^1}' +2 (\tau_0^2 - {\tau_0^3}')+ k_1^2 \tau_0^1+ 4 k_1 ({\tau_0^2}' - 3 k_1\tau_0^3),   
\end{align*}
 where we use \eqref{eq:ps03} and \eqref{eq:ps02}.
 Subtracting \eqref{eq:ps01} from \eqref{eq:ps12}, 
 \[
\tau_2^2 = 2 \tau_1^1 -\tau_0^0 +7 k_1 \tau_0^2 - 3 k_1 \tau_1^3 +2 \tau_0^3 + (\tau_1^2 - \tau_0^1)'.
 \]
 Then, by \eqref{eq:ps02} and \eqref{eq:ps03}, we can rephrase the 
 above equation as 
\begin{equation}\label{eq:X}
\tau_2^2 = 2 \tau_1^1 -\tau_0^0 +4 k_1 \tau_0^2 +2\tau_0^3 - 3 {k_1}' \tau_0^3 -6 k_1 
 {\tau_0^3}' + {\tau_0^2}''.
\end{equation}
Subtracting \eqref{eq:ps01} from \eqref{eq:ps23}, we have 
 \[
 \tau_1^1 + \tau_2^2 = \frac{1}{2}({\tau_0^1}- {\tau_2^3})'
 -2\tau_0^3 - 2 k_1 (\tau_1^3 + \tau_0^2).
 \]
 Then, by \eqref{eq:ps0123} and \eqref{eq:ps03}, we can 
 rephrase the above equation as 
\begin{equation}\label{eq:Y}
 \tau_1^1 + \tau_2^2 = \frac{1}{2}\left(-2{\tau_0^2}''- {\tau_0^3}'''\right)
 -2\tau_0^3 - 2 k_1 (2 \tau_0^2 + {\tau_0^3}').
\end{equation}
Then, subtracting \eqref{eq:Y}  from \eqref{eq:X}, we have 
\begin{equation}\label{eq:ps31100}
3 \tau_1^1 - \tau_0^0 = -8 k_1 \tau_0^2 -2 {\tau_0^2}'' + (3  { k_1}' -4)\tau_0^3 + 4k_1{\tau_0^3}' - \frac{1}{2} {\tau_0^3}'''.
\end{equation}
Now, adding \eqref{eq:ps11} and  \eqref{eq:ps22}, 
 \[
  \tau_3^2 - \tau_1^0  = 
 ( {\tau_1^1}' + {\tau_2^2}') + 3 k_1 (\tau_0^1 - \tau_2^3) -2  \tau_1^3 +2 \tau_0^2.
 \]
 By using \eqref{eq:Y} and \eqref{eq:ps0123}, we can rephrase 
 the above equation as 
\begin{equation}\label{eq:ps3210}
  \tau_3^2 - \tau_1^0  = 
-4 {k_1}' \tau_0^2 -10 k_1 {\tau_0^2}' - {\tau_0^2}''' -2 {\tau_0^3}'
 -2 {k_1}' {\tau_0^3 }' -5 k_1 {\tau_0^3}''
 - \frac{1}{2} {\tau_0^3}''''.
\end{equation}
Adding $4$\eqref{eq:ps01}, \eqref{eq:ps20} and \eqref{eq:ps31}, 
 we have 
\[
 8 \tau_0^0 = -(4 {\tau_0^1} + {\tau_2^0} + {\tau_3^1})' + 16 k_1 \tau_0^2 +8 
 \tau_0^3 + 4 k_1 (\tau_3^2 - \tau_1^0) - k_2 (\tau_0^1 - \tau_2^3).
\]
Then, by using \eqref{eq:ps3210} and \eqref{eq:ps0123}, we can rephrase 
 the above equation as
\begin{align}
 - 2 \tau_0^0 =  & \frac{1}{4}(4 {\tau_0^1} + {\tau_2^0} + {\tau_3^1}) '
 +4 k_1({k_1}' - 1) \tau_0^2 +
 \left(10 k_1^2 - \frac{1}{2} k_2\right) {\tau_0^2}'
 + k_1 {\tau_0^2}'''  \nonumber \\ &-2 \tau_0^3 
 +2 k_1\left(1 + {k_1}'\right) {\tau_0^3}' +\left( 5 k_1^2 - \frac{1}{4}k_2\right) 
 {\tau_0^3}'' + \frac{1}{2} k_1 {\tau_0^3}''''.
 \label{eq:ps002}
 \end{align}
 Adding \eqref{eq:ps002} and \eqref{eq:ps31100}, we have 
\begin{align}
 3 (\tau_1^1 - \tau_0^0) =& \frac{1}{4}(4 {\tau_0^1} + {\tau_2^0} + {\tau_3^1})' + 4 k_1
 ({k_1}' -3)\tau_0^2 + \left( 10 k_1^2 - \frac{1}{2}k_2\right) {\tau_0^2}'- 2 
 {\tau_0^2}''
 + k_1{\tau_0^2}'''  \nonumber \\ &+(3{k_1}' -6)\tau_0^3 
+2 k_1 (3  + {k_1}')  {\tau_0^3}'
+\left(5 k_1^2 -  \frac{1}{4}\right) {\tau_0^3}'' - \frac{1}{2} {\tau_0^3}'''+ \frac{1}{2} {k_1}'''' {\tau_0^3}.
\end{align}
 Therefore we have 
\begin{align}
 \tau_0^0 &- \tau_1^1 - 2\tau_0^3 - 4 k_1 \tau_0^2 \nonumber\\
= & - \frac{1}{12}(4 {\tau_0^1} + {\tau_2^0} + {\tau_3^1})' - \frac{4}{3} 
 k_1 {k_1}'\tau_0^2 - \frac{1}{6}( 20k_1^2 - k_2) {\tau_0^2}'+\frac{2}{3} 
 {\tau_0^2}'' -\frac{1}{3} k_1 {\tau_0^2}'''  \nonumber \\ 
 &- {k_1}' \tau_0^3 
-\frac{2}{3}k_1 ({k_1}' + 3) {\tau_0^3}'
- \frac{1}{12} \left(20 k_1^2 -  k_2 \right) {\tau_0^3}''
+ \frac{1}{6} {\tau_0^3}'''- \frac{1}{6} k_1 {\tau_0^3}''''.
 \label{eq:psintegrand}
\end{align}
 Finally, by using Stokes' theorem, we obtain
\begin{align*}
 \delta \int_C \omega  = 
\int_C  \left\{ 
- \frac{4}{3} k_1 {k_1}'\tau_0^2 - \frac{1}{6}( 20k_1^2 - k_2) {\tau_0^2}'  \right.  &-\frac{1}{3} k_1 {\tau_0^2}'''- {k_1}' \tau_0^3 -\frac{2}{3}k_1 (
{k_1}' + 3) {\tau_0^3}' \nonumber \\ 
 &  \left. 
- \frac{1}{12} \left(20 k_1^2 -  k_2 \right) {\tau_0^3}'' 
- \frac{1}{6} k_1 {\tau_0^3}''''
\right\} \om.
\end{align*}
 By using integration by parts, we finally obtain 
 differential equations for an extremal space curve as 
\begin{equation}\label{eq:extremalpspace}
\left\{
\begin{array}{l} 
\displaystyle 
{k_1}''' +16  k_1 {k_1}' -\frac{1}{2}{k_2}' =0, \\
\displaystyle 
{k_1}'''' +16 k_1 {k_1}'' +16 {({k_1}')}^2 + 6{k_1}' - \frac{1}{2} {k_2}''=0.
\end{array}
\right.
\end{equation}
 From these equations, the following theorem holds.
 \begin{theorem}[\cite{Ki}, pp. 233-234 in \cite{Bl}]
 \label{projspcextremal}
 A space curve without inflection points is projective extremal relative to 
 the length functional if and only if  the both curvatures are constant.
\end{theorem}
\noindent Proof.
 Subtracting the derivative of the first equation from the second equation \eqref{eq:extremalpspace}, we have ${k_1}' =0$. Then ${k_2}' =0$ follows immediately. 

\section{Classification of space curves with constant projective
curvatures} \label{appclass}

In Section \ref{subs-projinv}, we have seen that 
any nondegenerate
curve with constant projective curvatures is defined by
a differential equation
\begin{equation} \label{eqconst}
x''''= ax''+bx'+cx, 
\end{equation}
where the coefficients $a$, $b$ and $c$ are constants.
Such a curve is described by linearly independent solutions
and we get a classification of such curves by listing such solutions.

Let $e^{\lam t}$ be a solution of the differential equation, then
$\lam$ is a solution of the algebraic equation
\[\lam^4- a\lam^2 - b\lam - c=0.\]
Since the sum of the four roots is zero, we have the
following cases to consider separately.

\begin{tabular}{ll}
1. four distinct real roots: & $(\lam,\mu,\nu, -(\lam+\mu+\nu))$ \\
2. one set of double real roots: & $(2\lam, 2\mu, -\lam-\mu, -\lam-\mu)$ \\
3. one set of complex conjugate roots: & 
$(2\lam,2\mu,-(\lam+\mu)+ip,-(\lam+\mu)-ip)$ \\
4. two distinct sets of complex roots: & $(\lam+ip,\lam-ip,-\lam+iq, -\lam-iq)$ \\
5. one set of complex roots and double real roots:   & $(-\lam, -\lam, \lam+ip, \lam-ip)$ \\
6. two sets of double roots: & $( \lam, \lam, - \lam, - \lam)$ \\
7. pure imaginary roots: &  $(ip, ip, -ip, -ip)$ \\
8. triple roots: &  $(\lam, \lam, \lam, -3\lam)$ \\
9. trivial case: &  $(0,0,0,0)$,
\end{tabular}
\medskip

\noindent where $\lam$, $\mu$, $\nu$, $p$, $q$ are real constants.
For each case above, we compute the set of solutions to define
the immersion: we denote by CVi the curve in the case i above.
The associated differential equation defining the immersion is
also listed, 
where $a$, $b$, $c$  denote
the coefficients of the differential equation.
\medskip

{\small
\centering
\begin{tabular}{ll} \hline
{\small Curves} & mapping \\ \hline
CV1 & $[e^{-(\lam+\mu+\nu)t},e^{\lam t},e^{\mu t},e^{\nu t}]$ \\
CV2 & $[e^{-(\lam+\mu)t},te^{-(\lam+\mu)t},e^{2\lam t},e^{2\mu t}]$ \\
CV3 & $[e^{2\lam t},e^{2\mu t},e^{-(\lam+\mu)t}\cos(pt),e^{-(\lam+\mu)t}\sin(pt)]$ \\
CV4 & $[e^{\lam t}\cos(pt),e^{\lam t}\sin(pt), e^{-\lam t}\cos(qt),e^{-\lam t}\sin(qt)]$ \\
CV5 & $[e^{\lam t}\cos(pt),e^{\lam t}\sin(pt),e^{-\lam t},te^{-\lam t}]$ \\
CV6 & $[e^{\lam t}, t e^{\lam t}, e^{-\lam t}, t e^{-\lam t}]$ \\
CV7 & $[\cos(pt),\sin(pt),t\cos(pt),t\sin(pt)]$ \\
CV8 & $[e^{\lam t},te^{\lam t},t^2e^{\lam t},e^{-3\lam t}]$ \\
CV9 & $[1,t,t^2,t^3]$ \\
\hline
\end{tabular}
\smallskip

\begin{tabular}{llll} \hline
{\small Curves} & \qquad $a$ & \qquad $b$ & \qquad $c$ \\ \hline
CV1 & $\lam^2+\lam\mu+\lam\nu+\mu^2+\mu\nu+\nu^2$ & $-(\mu+\nu)(\nu+\lam)(\lam+\mu)$ & $\lam\mu\nu(\lam+\mu+\nu)$ \\
CV2 & $ 3\lam^2+2\lam\mu+3\mu^2$ & $2(\lam+\mu)(-\mu+\lam)^2$ & $-4\lam\mu(\lam+\mu)^2$ \\
CV3 & $ 3\lam^2+2\lam\mu+3\mu^2-p^2$ & $2(\lam+\mu)(\lam^2-2\lam\mu+\mu^2+p^2)$ 
  & $-4\lam\mu(\lam^2+2\lam\mu+\mu^2+p^2)$ \\
CV4 & $ 2\lam^2-q^2-p^2$ & $-2\lam(p-q)(p+q)$ & $-(\lam^2+q^2)(\lam^2+p^2)$ \\
CV5 & $2\lam^2-p^2$ & $-2p^2\lam$ & $-\lam^2(\lam^2+p^2)$ \\
CV6 & $ 2\lam^2$ & $0$ & $-\lam^4$ \\
CV7 & $ -2p^2$ & $0$ & $-p^4$ \\
CV8 & $ 6\lam^2$ & $-8\lam^3$ & $3\lam^4$ \\
CV9 & $ 0$ & $0$ & $0$ \\
\hline
\noalign{\medskip}
\multicolumn{4}{c}{Table A1:\   Space curves of constant projective curves}
\end{tabular}
}

\comment{
\begin{tabular}{lllll} \hline
{\small Curves} & mapping & \qquad A & \qquad B & \qquad C \\ \hline
CV1 & $[e^{-(a+b+c)t},e^{at},e^{bt},e^{ct}]$
 & $(a^2+ab+ac+b^2+bc+c^2)$ & $-(b+c)(c+a)(a+b)$ & $abc(a+b+c)$ \\
CV2 & $[e^{-(a+b)t},te^{-(a+b)t},e^{2at},e^{2bt}]$
 & $ (3a^2+2ab+3b^2)$ & $2(a+b)(-b+a)^2$ & $-4ab(a+b)^2$ \\
CV3 & $[e^{2at},e^{2bt},e^{-(a+b)t}\cos(pt),e^{-(a+b)t}\sin(pt)]$
 & $ (3a^2+2ab+3b^2-p^2)$ & $2(a+b)(a^2-2ab+b^2+p^2)$ 
  & $-4ab(a^2+2ab+b^2+p^2)$ \\
CV4 & $[e^{at}\cos(pt),e^{at}\sin(pt), e^{-at}\cos(qt),e^{-at}\sin(qt)]$
 & $ (2a^2-q^2-p^2)$ & $-2a(p-q)(p+q)$ & $-(a^2+q^2)(a^2+p^2)$ \\
CV5 & $[e^{at}\cos(pt),e^{at}\sin(pt),e^{-at},te^{-at}]$
 & $(2a^2-p^2)$ & $-2p^2a$ & $-a^2(a^2+p^2)$ \\
CV6 & $[\cos(pt),\sin(pt),t\cos(pt),t\sin(pt)]$
 & $ -2p^2$ & $0$ & $-p^4$ \\
CV7 & $[e^{at},te^{at},t^2e^{at},e^{-3at}]$
 & $ 6a^2$ & $-8a^3$ & $+3a^4$ \\
CV8 & $[1,t,t^2,t^3]$
 & $ 0$ & $0$ & $0$ \\
\hline
\end{tabular}
}

All curves, except CV4 and CV7, are
general-affine homogeneous and already were listed in
Section \ref{spcconst}. 

\newpage 

\subsection{$1$-parameter subgroups defining anharmonic curves}
 Thanks to Theorem \ref{spacenatural}, 
 each curve CVi is an orbit of a point $p$ under a $1$-parameter 
 subgroup $G$. We list them as follows:

{\small
\centering
\begin{tabular}{lcc}
Curves  & 1-parameter subgroup $G$  & Point $p$  \\ \hline
\noalign{\smallskip}
CV1 &\qquad  $\left(\begin{array}{cccc} e^{-(\lam+\mu+\nu)t} & & & \\
 & e^{\lam t}& & \\ & & e^{\mu t} & \\ & & & e^{\nu t}\end{array}\right)$ \qquad &
\qquad {\footnotesize $\left(\begin{array}{c} 1\\ 1 \\ 1 \\ 1 \end{array}\right)$} \\
CV2 &\qquad  $\left(\begin{array}{cccc} e^{-(\lam+\mu)t} &0 &0 &0 \\
te^{-(\lam+\mu)t} & e^{-(\lam+\mu)t}& 0  &0 \\ & & e^{2\lam t} & 0 \\ & &0 & e^{2\mu t}\end{array}\right)$ \qquad &
\qquad {\footnotesize $ \left(\begin{array}{c} 1\\ 0 \\ 1 \\ 1 \end{array}\right)$} \\
CV3 &\qquad  $\left(\begin{array}{cccc} e^{2\lam t} & & & \\
 & e^{2\mu t}& & \\ & & e^{-(\lam+\mu)t}\cos(pt) & e^{-(\lam+\mu)t}\sin(pt)\\ 
 & &e^{-(\lam+\mu)t}\sin(pt) & e^{-(\lam+\mu)t}\cos(pt)\end{array}\right)$ \qquad &
\qquad {\footnotesize $\left(\begin{array}{c} 1\\ 1 \\ 1 \\ 0 \end{array}\right)$ }\\
CV4 &\qquad  $\left(\begin{array}{cccc} e^{\lam t}\cos(pt) &-e^{\lam t}\sin(pt) & & \\
e^{\lam t}\sin(pt) & e^{\lam t}\cos(pt)& & \\ & & e^{-\lam t}\cos(qt) & e^{-\lam t}\sin(qt)\\ 
 & &e^{-\lam t}\sin(qt) & e^{-\lam t}\cos(qt)\end{array}\right)$ \qquad &
\qquad {\footnotesize $\left(\begin{array}{c} 1\\ 0 \\ 1 \\ 0 \end{array}\right)$ }\\
CV5 &\qquad  $\left(\begin{array}{cccc} e^{\lam t}\cos(pt) &-e^{\lam t}\sin(pt) & & \\
e^{\lam t}\sin(pt) & e^{\lam t}\cos(pt)& & \\ & & e^{-\lam t}& 0\\ 
 & & te^{-\lam t} & e^{-\lam t}\end{array}\right)$ \qquad &
\qquad{\footnotesize  $\left(\begin{array}{c} 1\\ 0 \\ 1 \\ 0 \end{array}\right)$ }\\
CV6 &\qquad  $\left(\begin{array}{cccc} e^{\lam t} &0  & & \\
t e^{\lam t} & e^{\lam t} & & \\ & & e^{-\lam t}& 0\\ 
 & & te^{-\lam t} & e^{-\lam t}\end{array}\right)$ \qquad &
\qquad {\footnotesize $\left(\begin{array}{c} 1\\ 0 \\ 1 \\ 0 \end{array}\right)$ }\\
CV7 &\qquad  $\left(\begin{array}{cccc} \cos(pt) &-\sin(pt) & & \\
\sin(pt) & \cos(pt)& & \\ t\cos(pt)& -t\sin(pt) & \cos(pt) & \sin(pt)\\ 
t\sin(pt) &t\cos(pt) &\sin(pt) & \cos(pt)\end{array}\right)$ \qquad &
\qquad {\footnotesize $\left(\begin{array}{c} 1\\ 0 \\ 0 \\ 0 \end{array}\right)$ }\\
CV8 &\qquad  $\left(\begin{array}{cccc} e^{\lam t} &0 &0 &0 \\
te^{\lam t} & e^{\lam t}& 0  &0 \\ {1\over 2}t^2 e^{\lam t} & te^{\lam t} & e^{\lam t}& 0 \\ 
0 & 0 &0 & e^{-3\lam t}\end{array}\right)$ \qquad &
\qquad {\footnotesize $\left(\begin{array}{c} 1\\ 0 \\ 0 \\ 1 \end{array}\right)$ }\\

CV9 &\qquad  $\left(\begin{array}{cccc} 1 &0 &0 &0 \\
t & 1 & 0  &0 \\{1\over 2}t^2 &t & 1 & 0 \\
{1\over 6}t^3 &{1\over 2}t^2 & t & 1\end{array}\right)$ \qquad &
\qquad {\footnotesize $\left(\begin{array}{c} 1\\ 0 \\ 0 \\ 0 \end{array}\right)$ }\\
\noalign{\smallskip}
\hline
\end{tabular}
 }

\vskip2pc

\noindent {\it Acknowledgment$:$} The authors are grateful to Professors
Masaaki Yoshida, Junichi Inoguchi, Wayne Rossman and Udo Hertrich-Jeromin 
for the helpful discussions 
given to us during the preparation of this manuscript.

\end{document}